\documentclass[reqno]{amsart}

\pdfoutput=1

\usepackage{amsfonts,amsmath,amsthm,amssymb}

\usepackage[square,numbers]{natbib}

\usepackage{graphicx}

\usepackage[colorlinks = true, citecolor = blue, urlcolor = blue]{hyperref}

\newtheorem{theorem}{Theorem}[section]
\newtheorem{lemma}[theorem]{Lemma}
\newtheorem{proposition}[theorem]{Proposition}
\newtheorem{corollary}[theorem]{Corollary}
\theoremstyle{definition}
\newtheorem{definition}[theorem]{Definition}
\newtheorem{remark}[theorem]{Remark}

\usepackage{geometry}
\usepackage{bbm}

\usepackage{tikz-cd}
\usepackage{mathtools}

\usepackage{algorithm}
\usepackage{algpseudocode} 

\makeatletter
\renewcommand{\varinjlim}{
\mathop{\mathpalette\varlim@{\rightarrowfill@\scriptscriptstyle}}\nmlimits@
}
\makeatother

\def \open {\underset{\clap{\scriptsize open}}{\subseteq}}
\def \closed {\underset{\clap{\scriptsize closed}}{\subseteq}}
\def \Cech {Čech\ }

\title[Sparse Eilenberg-MacLane Coordinates]{From Persistent Cohomology to Maps into Arbitrary Eilenberg-MacLane Spaces}

\author[Tony Xiaochen Xiao]{Tony Xiaochen Xiao}
\address{
\shortstack[l]{
Department of Mathematics, \\
College of Science, \\
Northeastern University, \\
Boston, MA, USA.}}
\email{xiao.xiaoc@northeastern.edu}

\author[Jose Perea]{Jose A. Perea}
\address{
\shortstack[l]{
Department of Mathematics, \\
College of Science, \\
Northeastern University, \\
Boston, MA, USA.}}
\email{j.pereabenitez@northeastern.edu}

\thanks{This work was partially supported by the National Science Foundation through award \# DMS-2603389 and CAREER award \# DMS-2415445}



\begin{document}

\bibliographystyle{abbrvnat}
\begin{abstract}
    In this paper, we address the construction of Eilenberg-MacLane coordinates — functions from datasets to Eilenberg-MacLane spaces $K(G,n)$ — that encode nontrivial persistent cohomology classes in arbitrary dimensions $n$ and coefficient group $G$. 
    The construction leverages the connecting homomorphism in the long exact sequence in \v{C}ech cohomology, from the short exact sequence of sheaves of continuous maps with values in $G$, $EG$ and $BG$.
    Crucially, the derivation is constructive: we give an explicit, recursive formula for preimages under the connecting homomorphism, which yields a constructive version of a classical cohomology vanishing theorem and, in turn, an algorithm for general Eilenberg-MacLane coordinates. 
    We validate our approach through experiments on synthetic datasets, and present novel examples of $K(\mathbb{Z}/2, 2)$-coordinates not addressed in prior work.
\end{abstract}

\maketitle
\section{Introduction}
\subsection{Background}
As a new field that has been developing rapidly since the 1990s, topological data analysis focuses on inferring topological features of an unknown topological space $\mathbb{X}$, given a finite sample $X=\{x_i\}_{i=1}^N \subseteq \mathbb{X}$. With the concept of persistent homology, introduced by Edelsbrunner, Letscher, and Zomorodian in 2002 \cite{edelsbrunner}, as well as the persistence barcodes (a visualizable invariant of persistent homology) and the efficient algorithm for computing it due to Carlsson and Zomorodian \cite{carlsson}, homological inference has become mature. Nowadays, given a finite   point cloud $X$ sampled from an underlying space $\mathbb{X}$ (such as $S^1, T^2, \mathbb{RP}^2$, etc.) embedded in a metric space $(\mathbb{M},d)$, the usual workflow is as follows:
\begin{enumerate}
\item  Build simplicial complexes, such as the Rips or \Cech complexes: 
\begin{equation}\label{def: Rips and Cech Complexes}
\begin{aligned}
\mathcal R_\alpha(X) &:= \{\sigma\subseteq X\mid \forall i,j\in \sigma,\ d(x_i,x_j)<\alpha\},\\
\check{\mathcal C}_\alpha(X) &:= \Big\{\sigma\subseteq X\ \Big|\ \bigcap_{x\in \sigma} B_{\alpha}(x)\neq \varnothing \Big\}.
\end{aligned}
\end{equation}
where $B_{\alpha}(x):=\{y\in \mathbb{M}\ |\ d(x,y)<\alpha\}$ denotes the ball with radius $\alpha$ centered at $x$.
Both constructions yield a filtered simplicial complex; that is, a 
family $\mathcal{K} = \{ \mathcal K_\alpha\}_{\alpha \in \mathbb{R}}$ of simplicial
complexes $\mathcal K_\alpha$ --- or $\mathcal K_j$, $j\in \mathbb{Z}$, if the index $\alpha$ is discretized to $\alpha_j$ --- so that $\mathcal K_\alpha \subseteq \mathcal K_{\alpha'}$ whenever $\alpha \leq \alpha'$.

\item Compute the persistent cohomology of $\mathcal{K}$, denoted   $PH^n(\mathcal{K};\mathbb{F}):=\bigoplus_{j\in \mathbb{Z}} H^n(\mathcal K_j;\mathbb{F})$. 
Since $X$ is finite,  $\bigcup \mathcal{K}$ will have finitely many simplices, and thus 
$PH^n(\mathcal{K}; \mathbb{F})$ can be endowed with the 
structure of a finitely generated graded module over the polynomial ring $\mathbb{F}[t]$, with coefficients in a field $\mathbb{F}$. By the classification theorem of finitely generated modules over PIDs,
$$PH^n(\mathcal{K};\mathbb{F})\cong \left(\bigoplus_{q=1}^Q t^{n_q}\cdot \mathbb{F}[t]\right)\oplus \left(\bigoplus_{l=1}^L (t^{m_l}\cdot \mathbb{F}[t])/(t^{m_l+d_l})\right)\label{classification of pid}$$
for unique $n_q, m_l\in \mathbb{Z}, d_l\in \mathbb{N}$. The multiset of \emph{persistence intervals} \[
[n_1,\infty), [n_2,\infty), ..., [n_Q,\infty), [m_1,m_1+d_1), [m_2,m_2+d_2), ..., [m_L,m_L+d_L)\]
is referred to as the barcode of $PH^n(\mathcal{K}; \mathbb{F})$, and provides a complete discrete invariant for the $\mathbb{F}[t]$-isomorphism type of $PH^n(\mathcal{K};\mathbb{F})$.
\item The long intervals in step 2 are in correspondence with generators $[\eta]\in H^n(\mathcal K_\alpha;\mathbb{F})$, and under appropriate hypotheses recover generators for the homology/cohomology of $\mathbb{X}$.
\end{enumerate}
\textit{Remark.} Typical software for computing persistence, such as Ripser~\cite{bauer2021ripser}, takes $\mathbb{F}:=\mathbb{Z}_q$ for a prime $q$ (throughout this paper, we write $\mathbb Z_q:=\mathbb Z/q\mathbb Z$, not to be confused with the $q$-adic integers).

\subsection{Problem} In many cases, the barcode of $PH^n(\mathcal K;\mathbb F)$ produced
by the workflow described above can be used reliably to infer the homology and cohomology of the underlying space $\mathbb X$ 
(this is known as \textit{homological/cohomological inference}~\cite[Chapter~4]{oudot2015persistence}).
Interpreting these classes, however, remains a central and challenging problem in data analysis: a barcode tells us that a class exists, not where it lives in the data. Concretely, given a persistent cohomology class $[\boldsymbol{\eta}]\in PH^n(\mathcal R(X);\mathbb Z_q)$, a natural question is: \textit{what information does $[\boldsymbol{\eta}]$ actually encode, and how do we understand it?} 
One  active direction  in addressing this question is to represent such a class by \emph{coordinates}, that is, maps from the data to a well-understood topological space. The philosophy is as follows:

At a suitable scale $\alpha$, the class $[\boldsymbol{\eta}]$ determines a class $[\eta]\in H^n(\mathcal K_\alpha;\mathbb Z_q)$, which can be interpreted as (or lifted to) a class $[\eta] \in H^n(\mathcal K_\alpha ; G)$ for appropriate $G$ finitely generated and Abelian. 
The representability of cohomology via Eilenberg-MacLane spaces, i.e.\ the natural bijection
$H^n(\mathcal K_\alpha; G) \cong [|\mathcal K_\alpha| , K(G,n)]$,
implies that $\eta$ is represented by a map $f_\eta : |\mathcal K_\alpha| \rightarrow K(G,n)$, unique up to homotopy. 
Moreover, if $X^{(\alpha)}:=\bigcup_{x\in X} B_{\alpha}(x)$ denotes the $\alpha$-neighborhood of $X$ in $\mathbb{M}$ (called the $\alpha$-offset/thickening of $X$), then when $\mathcal K_\alpha$ is either $\mathcal R_\alpha(X)$ or $\check{\mathcal C}_\alpha(X)$, one has a nerve map $X^{(\alpha)} \rightarrow |\mathcal K_\alpha|$ \cite{mccabe2026wassersteinstabilitynervemap} (which is a homotopy equivalence under appropriate conditions -- e.g., those of Leray's nerve lemma) that, after composing with $f_\eta$, yields a map $X^{(\alpha)} \rightarrow K(G,n)$.

This argument, however, is purely existential: $K(G,n)$ is in general an infinite-dimensional CW-complex without a canonical concrete model, and the bijection above does not by itself produce a computable map. Making $f_\eta$ explicit has so far been achieved only on particular cases of $K(G,n)$.

\subsection{Previous Work} These particular cases are as follows. In 2011, de Silva, Morozov, and Vejdemo-Johansson~\cite{desilva2011persistent} introduced circular coordinates (i.e.\ $X\to S^1 $) for data with nontrivial 1-dimensional persistent cohomology. 
In 2018, Perea~\cite{perea2017multiscale} presented a method for constructing maps to real and complex projective spaces   (i.e.\ $X^{(\alpha)}\to \mathbb{RP}^\infty, \mathbb{CP}^\infty  $). 
In 2019, Polanco and Perea~\cite{polanco2019coordinatizingdatalensspaces} developed lens coordinates (i.e., $X^{(\alpha)} \to S^\infty/\mathbb{Z}_q$) for cyclic coefficient groups $\mathbb{Z}_q$ and further developed lens-PCA. In 2020, Perea~\cite{perea2019sparse} worked out an explicit formula for sparse circular coordinates (i.e.\ $X^{(\alpha)}\to S^1$), in the sense that the map is also defined on an open neighborhood of the data. In 2023, Scoccola et al.~\cite{scoccola2023toroidal} identified a formal notion of geometric correlation between circular coordinates and described a systematic procedure for constructing low-energy torus-valued maps (i.e.\ $X^{(\alpha)}\to T^d$). In 2024, with spheres as targets, Schonsheck and Schonsheck~\cite{schonsheck2024spherical} constructed spherical coordinates $X\to S^2$ from second-degree persistent cohomology via variational optimization, and Nordwald, Garc\'ia-Redondo, and Monod~\cite{nordwald2026density} subsequently made this construction robust to sampling density.

The work above shares a common philosophy: in low degrees, $K(G,n)$ admits concrete geometric models --- e.g.\ $K(\mathbb Z,1)\simeq S^1$, $K(\mathbb Z_2,1)\simeq\mathbb{RP}^\infty$,
 $K(\mathbb{Z}_q , 1)\simeq L_q^\infty $, $K(\mathbb Z,2)\simeq\mathbb{CP}^\infty$ --- so that coordinates can be constructed directly into these space, and one can leverage case-specific optimizations. 
 This is why the past decade has proceeded case by case:

\begin{center}
\footnotesize
\begin{tabular}{|c|c|c|c|}
    \hline
    Coordinate Type & Cohomology Source & Target Space & Concrete Model\\
    \hline
    Circular \cite{desilva2009persistent, desilva2011persistent, perea2019sparse} & $PH^1(\mathcal{R}(X); \mathbb{Z})$ & $K(\mathbb{Z},1)$ & $S^1$\\
    Toroidal \cite{scoccola2023toroidal} & $PH^1(\mathcal{R}(X); \mathbb{Z}^d)$ & $K(\mathbb{Z}^d,1)$ & $T^d := S^1\times \cdots \times S^1$\\
    Real Projective \cite{perea2017multiscale} & $PH^1(\mathcal{R}(X); \mathbb{Z}_2)$ & $K(\mathbb{Z}_2,1)$ & $\mathbb{RP}^\infty$\\
    Lens \cite{polanco2019coordinatizingdatalensspaces} & $PH^1(\mathcal{R}(X); \mathbb{Z}_q)$ & $K(\mathbb{Z}_q,1)$ & $L_q^\infty := S^{\infty}/\mathbb{Z}_q$ \\
    Complex Projective \cite{perea2017multiscale, schonsheck2024spherical, nordwald2026density} & $PH^2(\mathcal{R}(X); \mathbb{Z})$ & $K(\mathbb{Z},2)$ & $\mathbb{CP}^\infty$ (or $S^2\subseteq \mathbb{CP}^\infty$) \\
    \hline
    Eilenberg-MacLane (\textbf{Our work}) & $PH^n(\mathcal{R}(X); G)$ & $K(G,n)$ & $B^nG$ ($n$-fold classifying space) \\
    \hline
\end{tabular}
\end{center}

\subsection{Our Contribution}
The purpose of this paper is to generalize these results to Eilenberg-MacLane coordinates $X^{(\alpha)}\to K(G,n)$, derived from persistent cohomology $PH^n(\mathcal R(X);G)$, for arbitrary dimensions $n\ge1$ and discrete Abelian groups $G$, including but not limited to $\mathbb{Z}$, $\mathbb{Z}_q$, and products thereof.
Here, by a persistent cohomology class 
$[\boldsymbol{\eta}] \in PH^n(\mathcal{K}; G)$ we 
mean a family $\{[\eta_\alpha]\}_{\alpha \in I}$
of cohomology classes $[\eta_\alpha] \in H^n(\mathcal{K}_\alpha ; G)$ indexed by an interval $I \subseteq \mathbb{R}$,
so that $[\eta_{\alpha'}]$ maps to 
$[\eta_\alpha]$, $\alpha \leq \alpha'$, under the homomorphism in cohomology induced by the inclusion $\mathcal{K}_\alpha \hookrightarrow \mathcal{K}_{\alpha'}$.

\noindent\textbf{1. A Unified Framework for Eilenberg-MacLane Coordinates.}

Our first contribution is methodological;  we propose to adopt a  model of $K(G,n)$ that applies uniformly: the iterated bar construction $B^nG$  of $G$ (Definition~\ref{Bar Construction}). 
The bar construction satisfies that if $G$ is an Abelian topological group, then so is $BG$. One may therefore iterate it to form $B^2G$, $B^3G$, etc. 
When $G$ has the discrete topology, one can show that $B^nG$ is a $K(G,n)$ (Proposition~\ref{A good model for K(G,n)}).

Second, given a cohomology class in degree $n$ and coefficients in $G$, we use the long exact sequence in \v{C}ech cohomology induced by the universal bundle $B^{k-1}G\to E(B^{k-1}G)\to B^k G$, $k = 1,\ldots, n$, to lower the degree of cohomology and increase the bar iterates (Section~\ref{sec:lowering-degree-general}), one connecting homomorphism at a time, down to cohomological degree one or zero.
What makes this possible is the following theorem. 
For $B$ paracompact Hausdorff and $E$ a contractible topological Abelian group, we have\footnote{Here $\mathcal F_E$ denotes the sheaf of continuous $E$-valued maps on $B$, as defined in Section~\ref{sec: Soft Sheaves and the Vanishing of Cech Cohomology}.} $\check H^n(B;\mathcal F_E)=0$ for $n>0$.
This vanishing result is classical~\cite[Lemma~3.2 and \S4]{huber1961homotopical},
but it contains the core algorithmic difficulty we address in this work:
in order to obtain an implementation of the above degree lowering strategy, one would require an algorithm that, from a degree-$n$ \v{C}ech cocycle with coefficients in $\mathcal{F}_E$, explicitly constructs a degree-$(n-1)$ \v{C}ech cochain whose coboundary is that cocycle.
Theorem~\ref{thm:intro-Qn} gives, to our knowledge, the first constructive proof, and moreover one at the level of open covers.

\begin{theorem}[see Theorem~\ref{thm:Qn-correctness}]\label{thm:intro-Qn}
Let $B$ be a metric space with a finite cover $\mathcal U$ by open balls, and $E$ a contractible topological Abelian group. The $Q^n$-Algorithm (Algorithm~\ref{alg: Q^n-Algorithm}) 
$$Q^n:\check Z^n(\mathcal U;\mathcal F_E)\to\check C^{n-1}(\mathcal U;\mathcal F_E)$$
assigns to every \v{C}ech-cocycle $\eta$ an explicit cochain $Q^n(\eta)$ with $d^{n-1}(Q^n(\eta))=\eta$. In particular, $\check H^n(\mathcal U;\mathcal F_E)=0$ for all $n>0$.
\end{theorem} 

\noindent\textbf{2. Design of the Sparse Eilenberg-MacLane Algorithm.}
As a corollary of Theorem~\ref{thm:intro-Qn}, we obtain the \emph{Sparse Eilenberg-MacLane Algorithm}:

\begin{corollary}[see Section~\ref{sec: Eilenberg-MacLane Algorithm}]\label{cor:intro-EM}
Let $  K$ be a finite simplicial complex, and $G$ a discrete Abelian group. 
For every simplicial cocycle $\eta\in Z^n( K; G)$, the Eilenberg-MacLane Algorithm (Algorithm~\ref{alg: Eilenberg-MacLane Algorithm}), applied to the cover of $| K|$ by open vertex stars, outputs an explicit map $f_\eta: |  K|\to K(G,n)$ whose homotopy class corresponds to $[\eta]$ under the isomorphism $H^n(  K; G)\cong[|  K|,K(G,n)]$.
\end{corollary}

In the setting of persistent cohomology with coefficients in  $\mathbb{Z}_q$ for $q\geq 2$ prime, let $[\boldsymbol{\eta}]\in PH^n(\mathcal R(X);\mathbb Z_q)$ be an $n$-dimensional persistent cohomology class for the Rips filtration on $X$, detected with persistence interval $[a,b)$ satisfying $a<\frac b2$. Then, for any choice of $\alpha\in[a,\frac b2)$, the algorithm  outputs a map $X^{(\alpha)}\to B^nG$ into the model $B^nG$ of $K(G,n)$, where the coefficient group $G$ is determined by Bockstein lifting (Section~\ref{sec:bockstein-lifting}). The algorithm has the following key features:

\emph{Explicit formulas:} Each step of the construction evaluates  a explicit contracting homotopy (Equation \ref{eq:contracting-homotopy}) and piecewise-linear bump functions on the cover; no  optimization, or global solve is required. 
The only linear-algebraic step in the pipeline is the Bockstein lifting, which in most cases of interest (e.g., $G$ finitely generated) reduces to integer linear systems solvable via Smith normal form.

\emph{Sparse structure:} Our coordinate is defined not merely on the data $X$ itself, but on the entire open neighborhood $X^{(\alpha)}$ (thickening of $X$).

\emph{Milnor variant:} We also designed a variant of this pipeline, the \emph{Milnor variant} (Section~\ref{sec: Milnor-variant}), which replaces the last bar construction by the Milnor construction $\mathcal{B}(\cdot)$ and instead outputs a map $X^{(\alpha)}\to\mathcal B(B^{n-1}G)$, into another model of $K(G,n)$. It recovers the circular coordinates of~\cite{perea2019sparse} for $n=1$, and, whenever $B^{n-1}G$ can be replaced by a metric group, one can endow the output coordinate space with an explicit metric, making distance-based methods such as Isomap directly applicable.

\noindent \textbf{3. Experimental Validation and Novel Examples.}
We validate our framework through five experiments, organized into three groups, each isolating a distinct property of the construction:

\textit{Coefficient choice ($X\approx S^1$).} We compute both a $K(\mathbb Z,1)$- and a $K(\mathbb Z_q,1)$-coordinate from the same persistent cohomology class, reproducing the sparse  circular coordinate construction \cite{perea2019sparse} as the special case where the canonical pipeline and the Milnor variant coincide. 
This comparison demonstrates concretely why Bockstein lifting matters: the two coordinates are not interchangeable, and choosing the wrong coefficient group yields a needlessly more complex, less informative representation of the same underlying data.

\textit{Higher-degree iteration and its cost ($X\approx S^2$).} We compute a $B^2\mathbb Z$-coordinate via the canonical pipeline — confirming that, unlike the classical principal-bundle approach, our construction genuinely iterates to arbitrary degree $n$ — and a $\mathcal B(B\mathbb Z)$-coordinate via the Milnor variant from the same class. Both exhibit a geometric distortion intrinsic to the construction, which we show can be corrected after the fact via Isomap, restoring the expected persistent cohomology of $S^2$.

\textit{A genuinely new coordinate type ($X\approx M(\mathbb Z_q,1)$).} We construct a $\mathcal B(B\mathbb Z_q)$-coordinate from a torsion class in $H^2(M(\mathbb Z_q,1);\mathbb Z_q)$ — a degree-two, $\mathbb Z_q$-coefficient coordinate for a Moore space that is neither a lens coordinate \cite{polanco2019coordinatizingdatalensspaces} nor a projective coordinate \cite{perea2017multiscale}, and has not appeared in prior work. 

Together, these experiments show that our method (i) subsumes the classical, degree-one constructions as a special case, (ii) is the first to iterate to arbitrary degree $n$ in practice, and (iii) computes coordinate types that no prior framework computes — while also exposing, and correcting, a distortion intrinsic to the construction itself.

\subsection{An Overview of the Sparse Eilenberg-MacLane Algorithm}\label{sec: sparse-eilenberg-maclane-algorithm}

The design of the Sparse Eilenberg-MacLane Algorithm is developed in detail in Section~\ref{sec: methodology}. 
So that readers can understand the main results without the sheaf-theoretic machinery, we present the algorithm here in the language of functions and \v{C}ech cocycles alone.

\begin{enumerate}

\item Let $X$ be the input data set, a finite metric space with metric $d$.\footnote{One may optionally work with a subset of landmarks $L=\{\ell_1,\ldots,\ell_N\}\subset X$ instead of $X$ itself, e.g.\ chosen via \texttt{maxmin} sampling, so as to reduce the size of the input; in what follows we simply write $X$ for whichever set is used. In this case one additionally requires $\alpha>r_L:=\max_{x\in X}\min_{\ell\in L}d(x,\ell)$, so that $L^{(\alpha)}$ contains all data points.} 
Suppose that an $n$-dimensional persistent cohomology class $[\boldsymbol{\eta}]\in PH^n(\mathcal R(X);\mathbb Z_q)$ for the Rips filtration on $X$, with coefficients in $\mathbb Z_q$ for a prime $q$, is detected with persistence interval $[a,b)$ satisfying $a<\frac b2$. 
Choose $a\leq \alpha<\frac b2$ and let $\eta'\in Z^n(\mathcal R_{2\alpha}(X);\mathbb Z_q)$ be a cocycle representative for $[\boldsymbol{\eta}]$.

\item Lift $\eta' \in Z^n(\mathcal R_{2\alpha}(X);\mathbb{Z}_q)$ to a cocycle
    $\eta'' \in Z^n(\mathcal R_{2\alpha}(X);G)$. That is, one for which $\kappa_*([\eta'']) = [\eta']$ in the long exact sequence:
    \[\cdots\to H^n(\mathcal R_{2\alpha}(X); H)\xrightarrow{\iota_*} H^n(\mathcal R_{2\alpha}(X); G)\xrightarrow{\kappa_*} H^n(\mathcal R_{2\alpha}(X); \mathbb{Z}_q)\xrightarrow{\Delta} H^{n+1}(\mathcal R_{2\alpha}(X); H)\to \cdots\]
induced by a choice of short sequence $0 \to H \to G \to \mathbb{Z}_q\to 0$.
To this end, compute the Bockstein algorithm to determine if $\Delta([\eta'])=0$ for the sequences in Table~\ref{tab:bockstein-tower} (Section~\ref{sec:bockstein-lifting}), and take $G$ from the last row that returns ``success''. (In practice, one first tests the row $\infty$; if it succeeds, take $G=\mathbb{Z}$. Otherwise, test the rows $0,1,2,\dots$ in turn until one fails; the row immediately preceding it is the last success.)

Each test $\Delta([\eta'])=0$ amounts to deciding whether a specific integer cochain lies in the image of a coboundary matrix, and hence is solvable via Smith normal form.

\item Pass from $\eta''\in Z^n(\mathcal R_{2\alpha}(X);G)$ to $\eta \in \check Z^n(\mathcal U^\alpha;\mathcal F_G)$, where $\mathcal U^\alpha:=\{U_i:=B_\alpha(x_i)\}$ (see Section~\ref{sec:from-simplicial-to-cech} for details): for each $n$-simplex $i_0<\cdots<i_n$, let $\eta_{i_0\cdots i_n}:U_{i_0\cdots i_n}\to G$ be the  map which is constant and equal to $\eta''_\sigma$, where $\sigma := \{x_{i_0},\ldots,x_{i_n}\}$. (This is well defined, since $U_{i_0\cdots i_n}\neq\varnothing$ implies $\sigma\in\check{\mathcal C}_\alpha(X)\subseteq\mathcal R_{2\alpha}(X)$.)

\item Fix an open cover filtration of $X^{(\alpha)}$: $\mathcal U^\alpha=\mathcal U^{\gamma_0}\to\mathcal U^{\gamma_1}\to\cdots\to\mathcal U^{\gamma_{N-1}}$, where $N$ is the number of open sets in $\mathcal U$. The easiest way is to choose a descending sequence of radii $\alpha=\gamma_0>\cdots>\gamma_{N-1}$ and define $\mathcal{U}^r:=\{B_r(x)\}_{x\in X}$, chosen so that $\mathcal U^{\gamma_{N-1}}$ still covers $X^{(\alpha)}$. Set up bump functions $g_1,g_2,\ldots,g_{N-1}:X^{(\alpha)}\to[0,1]$ (in Section~\ref{sec: the-Q^n-Algorithm}) and the contracting homotopy $h:EG\times[0,1]\to EG$ of $EG$ (in Proposition~\ref{Contractibility of EG}).

\item Set $H:=G$, $k:=n$, and $\eta^{(n)}:=\eta$. While $k\geq 1$, repeat:
\begin{enumerate}
    \item Replace each component $\eta^{(k)}_{i_0\cdots i_k}:U_{i_0\cdots i_k}\to H$ by the $EH$-valued map $x\mapsto\bigl(\eta^{(k)}_{i_0\cdots i_k}(x),1\bigr)$, i.e., a $0$-simplex of $EH$ (Definition~\ref{Bar Construction}); components with a repeated index are taken to be the identity $e_{EH}$.
    \item For each $(k-1)$-simplex $i_0<\cdots<i_{k-1}$, set $t^{(0)}_{i_0\cdots i_{k-1}}:=\eta^{(k)}_{0i_0\cdots i_{k-1}}$ on $U_0\cap U_{i_0\cdots i_{k-1}}$. Then, for $i=1,2,\dots,N-1$, define $t^{(i)}_{i_0\cdots i_{k-1}}:\tilde U_i^{\gamma_i}\cap U_{i_0\cdots i_{k-1}}\to EH$ by

    \fbox{
    $\footnotesize t^{(i)}_{i_0\cdots i_{k-1}}(x) := \left\{ \begin{array}{ll}
        t^{(i-1)}_{i_0\cdots i_{k-1}}(x) & \text{if } x\in \tilde{U}_{i-1}^{\gamma_i}\cap U_{i_0\cdots i_{k-1}} \\
        \eta^{(k)}_{ii_0\cdots i_{k-1}}(x) + \sum\limits_{l=0}^{k-1} (-1)^l\cdot h\bigl(t^{(i-1)}_{ii_0\cdots \hat{i_l}\cdots i_{k-1}}(x), g_i(x)\bigr) & \text{if } x\in U_i^{\gamma_i}\cap U_{i_0\cdots i_{k-1}}
    \end{array}\right.$
    }

    where $\tilde{U}_i^r:=U_0^r\cup U_1^r\cup \cdots \cup U_i^r$.
    \item Set $\eta^{(k-1)} := \bigl( j\circ t^{(N-1)}_{i_0\cdots i_{k-1}} \bigr)_{i_0<\cdots<i_{k-1}}$, where $j:EH\to BH$ is the quotient map (Definition~\ref{Bar Construction}). Since $\tilde U_{N-1}^{\gamma_{N-1}}=X^{(\alpha)}$, each $t^{(N-1)}_{i_0\cdots i_{k-1}}$ is defined on all of $U_{i_0\cdots i_{k-1}}$, and $\eta^{(k-1)}$ defines a \Cech $(k-1)$-cocycle with values in $BH$.
    \item Set $H:=BH$ and $k:=k-1$.
\end{enumerate}

\item \footnote{As an alternative (the \emph{Milnor variant}), using a partition of unity $\{\phi_j\}$ subordinate to $\mathcal U^\alpha$, one may apply formula~\ref{eq:milnor-formula} to the \Cech $1$-cocycle $\eta^{(1)}$; this yields a map $X^{(\alpha)}\to\mathcal B(B^{n-1}G)$ into Milnor's classifying space of $B^{n-1}G$, which is again a $K(G,n)$, being homotopy equivalent to $B^nG$.} As a \Cech $0$-cocycle, $\eta^{(0)}$ (where $\eta^{(0)}_{i_0}: U_{i_0}\to B^nG$) can be naturally glued to a map $f_\eta:X^{(\alpha)}\to B^nG$. Since $B^nG$ is a $K(G,n)$ (Proposition~\ref{A good model for K(G,n)}), we get a map $f_\eta:X^{(\alpha)}\to K(G,n)$.

\renewcommand{\algorithmicrequire}{\textbf{Input:}}
\renewcommand{\algorithmicensure}{\textbf{Output:}}

\end{enumerate}

\subsection{Organization} Section~2 develops the preliminaries needed for the rest of the paper: an explicit model for $K(G,n)$ via the iterated bar construction, and the long exact sequence of \Cech cohomology underlying the degree-lowering mechanism central to this work. Section~3 develops the full pipeline transporting a persistent cohomology class $[\eta]\in PH^n(\mathcal R(X);\mathbb Z_q)$ into an explicit coordinate $X^{(\alpha)}\to K(G,n)$, together with an alternative route (the Milnor variant) through the classical theory of principal bundles. Section~4 makes this construction algorithmic, introducing the $Q^n$-Algorithm and assembling it, together with its Milnor variant, into the complete Eilenberg-MacLane Algorithm. Section~5 presents five experiments validating and illustrating this pipeline on $S^1$, $S^2$, and a Moore space $M(\mathbb Z_q,1)$, including a comparison with the classical, principal-bundle-based construction of circular coordinates, and a novel $\mathcal B(B\mathbb Z_q)$-coordinate not treated in prior work. Section~6 concludes with directions for future work.

\section{Preliminaries}

\subsection{A good model for $K(G,n)$.} \label{section: a good model for K(G,n)}

\begin{definition}[Eilenberg-MacLane Spaces]
\label{Eilenberg-MacLane Spaces}
Let $G$ be a group and let $n\geq 1$ be an integer. 
An \textit{Eilenberg-MacLane space} $K(G,n)$ is a connected CW-complex whose homotopy groups are all trivial, except in degree $n$, where it is isomorphic to $G$. 
That is, $\pi_i(K(G,n))=\left\{\begin{array}{ll}
    G & (\text{if } i=n) \\
    0 & (\text{else})
\end{array}\right.$
\end{definition}

Given $G$ and $n$, there are numerous concrete models for constructing $K(G,n)$, including the iterated bar construction~\cite[\S4]{milgram1967bar}, the Dold--Thom construction via infinite symmetric products of Moore spaces~\cite{dold1958quasifaserungen}, \cite[\S4.K]{hatcher2002algebraic}, or the direct approach of attaching cells to a wedge sum of copies of $S^n$ to successively kill the higher homotopy groups~\cite[\S4.2]{hatcher2002algebraic}, among others. Here we take the iterated bar construction route  — $BG, B^2G, ..., B^{n-1}G$ —  and present its explicit details in the following paragraphs.

\begin{definition}[Classifying Spaces]
\label{Classifying Spaces} For a topological group $G$, let $\underline{E}G$ be a weakly contractible space  with a free continuous action $\underline{E}G\times G\to \underline{E}G$. 
Then, the quotient $\underline{B}G:=\underline{E}G/G$  yields a principal $G$-bundle $G\to \underline{E}G\to \underline{B}G$, called the \textit{universal bundle} of $G$, and  
the space $\underline{B}G$ is called a \textit{classifying space for $G$}.
\end{definition}

Although the construction of $\underline{E}G$ is not unique, it turns out that all classifying spaces for   $G$ are homotopy equivalent to each other.

    \begin{definition}[Bar Construction $BG$]
    \label{Bar Construction} 
    Let $\Delta_n := \{(t_0,...,t_n) \in \mathbb{R}^{n+1} \ | \ t_i\geq 0, t_0+...+t_n=1\}$ denote the standard $n$-simplex. For a topological group $G$ with identity $e\in G$, let
        $$EG:=\bigsqcup_{n\geq 0} G^{n+1}\times \Delta_n/\sim$$
    with the equivalence relation $\sim$ given by
\begin{align*}
((g_0,\dots,g_n),(t_0,\dots,t_n))
&\sim ((g_1,\dots,g_n),(t_1,\dots,t_n)) && \text{if } t_0=0,\\
&\sim ((\dots,g_{i-2},g_{i-1}g_i,g_{i+1},\dots),(\dots,t_{i-1},t_{i+1},\dots)) && \text{if } t_i=0,\ i>0,\\
&\sim ((\dots,g_{i-1},g_{i+1},\dots),(\dots,t_{i-1},t_i+t_{i+1},t_{i+2},\dots)) && \text{if } g_i=e,\ i<n.
\end{align*}
where the first two relations are the face identifications and the last is the degeneracy identification.
    Here, each $G^{n+1}\times\Delta_n$ is endowed with the product topology, the disjoint union $\bigsqcup_{n\geq 0}G^{n+1}\times\Delta_n$ carries the coproduct topology, and $EG$ is equipped with the resulting quotient topology under $\sim$.

The right $G$-action on $EG$ is given by: $$\begin{array}{ccccc}
    EG & \times & G & \to & EG \\
    \ [(g_0,g_1,...,g_n),(t_0,t_1,...,t_n)] & \cdot & g & := & [(g_0,g_1,...,g_ng),(t_0,t_1,...,t_n)]
\end{array}$$
Since $EG$ is weakly contractible (Proposition~\ref{Contractibility of EG}), then the orbit space $BG:=EG/G$  is indeed a classifying space for $G$, referred to as the \textit{bar construction} of $G$.
    \end{definition}

\begin{definition}[Group law on $BG$]
    \label{group law on BG}
Regard a point in 
$BG$ as   the $G$-orbit of an equivalence class $[(g_0,...,g_n),(t_0,...,t_n)] \in EG$, with representative barycentric coordinates $(t_0,t_1, \ldots, t_n) \in \Delta_n$. 
We equip $BG$ with a group structure by defining, for two representative points with the same barycentric coordinates, as the $G$-orbit of 
$$[(g_0,...,g_n),(t_0,...,t_n)]\cdot[(g'_0,...,g'_n),(t_0,...,t_n)] := [(g_0g'
_0,...,g_ng'_n),(t_0,...,t_n)]$$
where the group operation is applied coordinate-wise in $G$. 
Given two points in $EG$ one may always find representatives with the same barycentric coordinates — obtained by repeatedly applying the degeneracy identification to introduce copies of $e$ at the appropriate coordinates (see the proof of \cite[Theorem 14.3]{petermaySimplicial}) — at which point the coordinate-wise multiplication above applies.
\end{definition}

\textit{Remark 1.} One may check that this multiplication is well-defined on equivalence classes and respects the relations defining $EG$, and is therefore compatible with the identifications used to form $BG$.

\textit{Remark 2.} It is also important to point out that this group law preserves the Abelian-ness, i.e., if $G$ is Abelian, then so is the bar construction $BG$ — this is immediate from Definition~\ref{group law on BG}, since the multiplication is applied coordinate-wise within $G$ itself.\\

As a result, under the condition that $G$ is an Abelian topological group, one may iteratively build $n$-fold bar constructions: $BG$, $B^2G$, $B^3G$, ..., etc., and they are all Abelian topological groups.

\begin{proposition}[Contractibility of $EG$]
\label{Contractibility of EG}
The space $EG$ of Definition~\ref{Bar Construction} is contractible. Explicitly, a contracting homotopy $h: EG\times [0,1]\to EG$ is given by
\begin{equation}\label{eq:contracting-homotopy}
\begin{array}{rcl}
    h: \left( \bigsqcup_{n\geq 0} G^{n+1}\times \Delta^n / \sim \right) \times [0,1] & \to & \bigsqcup_{n\geq 0} G^{n+1}\times \Delta^n / \sim \\
    ((g_0,\dots,g_n),(t_0,\dots,t_n),t) & \mapsto & (((g_0\cdots g_n)^{-1},g_0,\dots,g_n),(1-t,tt_0,\dots,tt_n)),
\end{array}
\end{equation}
which satisfies $h(\cdot,1)=\mathrm{id}$ and $h(\cdot,0)\equiv e$.
\end{proposition}

\begin{proof}
The contractibility of $EG$ is classical (cf.~\cite[Chapter~16, \S5]{petermay}), but we provide a specific null-homotopy as it will be needed in the algorithms later on. 
That $h$ is a well-defined continuous contracting homotopy is verified in Appendix~\ref{app:contracting-homotopy}.
\end{proof}

We finish this section by indicating why the $n$-fold bar construction is indeed a model for $K(G,n)$.

\begin{proposition}[A good model for $K(G,n)$; {\cite[\S4]{milgram1967bar}}]
\label{A good model for K(G,n)}
If $G$ is a discrete Abelian topological group, then for any $n\geq 1$, $B^nG$ is a $K(G,n)$.
\end{proposition}

\begin{proof}
For any topological Abelian group $G$, the long exact sequence of homotopy groups of the fibration $G\to EG\to BG$, with $EG$ contractible (Proposition~\ref{Contractibility of EG}), gives isomorphisms $\pi_{i+1}(BG)\cong\pi_i(G)$ for all $i\geq0$, and $BG$ is path-connected as a quotient of $EG$. Applying this successively to $G, BG,\dots,B^{n-1}G$, with $G$ discrete, yields $\pi_n(B^nG)\cong G$ and $\pi_i(B^nG)=0$ for $i\neq n$.
\end{proof}

\subsection{\Cech Resolutions} \label{sec: Cech Resolutions}

Before turning to the principles and structure of the Eilenberg-MacLane Algorithm itself, we must first revisit Čech cohomology in a more general form. Readers familiar with Čech cohomology will most likely recognize it in the following form: given an open cover $\mathcal{U}$ of a space $B$ and a fixed sheaf $\mathcal{F}$ of abelian groups on $B$, the cochain groups $\check{C}^n(\mathcal{U},\mathcal{F})$, equipped with the usual coboundary maps, form a cochain complex whose cohomology is $\check{H}^n(\mathcal{U},\mathcal{F})$. In this picture, each $\check{C}^n(\mathcal{U},\mathcal{F})$ is treated simply as an abelian group attached to the fixed space $B$. 
For the arguments developed in  Section \ref{sec: methodology}    we will need to speak of sections of $\check{C}^n(\mathcal{U},\mathcal{F})$ over arbitrary open and closed subsets of $B$ — in particular, to invoke sheaf-theoretic properties such as softness, which are meaningful only for sheaves and have no counterpart for a group attached to a single fixed space. 
We   begin  by recasting each cochain group as a sheaf in its own right.

\begin{definition}[\Cech Sheaves]
    Let $B$ be a topological space and $\mathcal{U}=\{U_i\}_{i\in I}$   an open cover of $B$. 
    Let $\mathcal{F}$ be a sheaf of Abelian groups over $B$ and $n$  a non-negative integer. 
    Then, the assignment
    $$V \mapsto \check{C}^n(\mathcal{U},\mathcal{F})(V):= \prod_{i_0<...<i_n} \mathcal{F}(U_{i_0...i_n}\cap V),\qquad \forall V\open B$$
    defines a sheaf over $B$, together with a sheaf map (i.e. sheaf morphism)
    $d^n: \check{C}^n(\mathcal{U},\mathcal{F})\to \check{C}^{n+1}(\mathcal{U},\mathcal{F})$ given by
    $$\begin{array}{rcl}
    d^n: \check{C}^n(\mathcal{U},\mathcal{F})(V) & \to & \check{C}^{n+1}(\mathcal{U},\mathcal{F})(V) ,\qquad \forall V\open B\\
    (s_{i_0...i_n})_{i_0<...<i_n} & \mapsto & \left( \sum_{j=0}^{n+1} (-1)^j s_{i_0...\hat{i_j}...i_{n+1}}\vert_{U_{i_0...i_{n+1}}\cap V} \right)_{i_0<...<i_{n+1}} \\
    \end{array}$$
    The sheaf $\check{C}^n(\mathcal{U},\mathcal{F})$ is called the \textit{\Cech sheaf} of $\mathcal{U}$ with coefficients in $\mathcal{F}$ and degree $n$.
\end{definition}

\begin{lemma}[\Cech Resolution {\cite[Lemma~III.4.2]{hartshorne1977algebraic}}]
\label{Cech Resolution}
The maps $d^n$ form an exact sequence of sheaves
\[0\to \mathcal{F}\xrightarrow{e} \check{C}^0(\mathcal{U},\mathcal{F})\xrightarrow{d^0} \check{C}^1(\mathcal{U},\mathcal{F})\xrightarrow{d^1} \cdots
\xrightarrow{d^{n-1}}\check{C}^n(\mathcal{U},\mathcal{F})\xrightarrow{d^n}\cdots,\]
called the \textit{\Cech resolution} of $\mathcal{F}$.
\end{lemma}

\begin{definition}[\Cech Cohomology]
    Let $B$ be a topological space,  $\mathcal{U}$   an open cover of $B$, and $\mathcal{F}$ a sheaf of Abelian groups over $B$. 
    The \textit{\Cech Cohomology groups} of $B$ over $\mathcal{U} $ and coefficients in $\mathcal{F}$ are defined as the cohomology of the cochain complex formed by taking global sections of the \Cech Resolution, i.e.,
    $$\check{H}^n(\mathcal{U},\mathcal{F}) := H_n(\mathcal{C^\bullet}):=\ker d^n_B/\text{im } d^{n-1}_B, \qquad \forall n\geq 1$$
    where
    $$\mathcal{C^\bullet}: 0\to \mathcal{F}(B)\xrightarrow{e_B} \check{C}^0(\mathcal{U},\mathcal{F})(B)\xrightarrow{d^0_B} \check{C}^1(\mathcal{U},\mathcal{F})(B)\xrightarrow{d^1_B} \cdots  
    \xrightarrow{d^{n-1}_B}\check{C}^n(\mathcal{U},\mathcal{F})(B) \xrightarrow{d^n_B} \cdots$$
\end{definition}

\textit{Remark.} By convention, $\check{H}^0(\mathcal{U},\mathcal{F}):=\ker d^0_B$; moreover, since global sections is a left-exact functor, Lemma~\ref{Cech Resolution} implies $\ker d^0_B=\mathrm{im}\ e_B$, so that $e_B$ induces an isomorphism $\check{H}^0(\mathcal{U},\mathcal{F})\cong\mathcal{F}(B)$.\\

Before stating the result, it is worth recalling that Čech cohomology, as defined above, depends a priori on the choice of open cover $\mathcal{U}$ of $B$ and not merely on $B$ itself; different covers may in general yield different groups $\check{H}^n(\mathcal{U},\mathcal{F})$. However, one may approach the cohomology intrinsic to $B$ by successively refining $\mathcal{U}$: in general, one defines
$$\check{H}^n(B,\mathcal{F}) := \varinjlim_{\mathcal U} \check{H}^n(\mathcal{U},\mathcal{F})$$
the direct limit taken over all open covers $\mathcal{U}$ of $B$, ordered by refinement. The following lemma addresses the coarsest possible case in this limit — namely, when $B$ itself already appears as a member of $\mathcal{U}$ — in which case the resulting Čech cohomology vanishes in every positive dimension.

\begin{lemma}[Local fineness of \Cech Cohomology]
\label{Local fineness of Cech Cohomology}
If $B\in\mathcal U$, then $\check H^n(\mathcal U,\mathcal F)=0$ for all $n\geq1$.
\end{lemma}

\begin{proof}
Let $B=U_i$. For a cocycle $\eta\in\check C^n(\mathcal U,\mathcal F)(B)$, the cochain $\tau\in\check C^{n-1}(\mathcal U,\mathcal F)(B)$ given by $\tau_{i_0\cdots i_{n-1}}:=\eta_{i\,i_0\cdots i_{n-1}}$ (well defined since $U_{i\,i_0\cdots i_{n-1}}=U_{i_0\cdots i_{n-1}}$, and using the usual alternating convention for repeated or unordered indices) satisfies $d^{n-1}_B\tau=\eta$ by the cocycle condition $d^n_B\eta=0$; see the proof of \cite[Lemma~III.4.2]{hartshorne1977algebraic}.
\end{proof}

\subsection{Soft Sheaves and the Vanishing of \Cech Cohomology} \label{sec: Soft Sheaves and the Vanishing of Cech Cohomology}

Having set up \Cech resolutions and \Cech cohomology in Section~\ref{sec: Cech Resolutions}, we now prove the following vanishing result. 
The result itself is known, but to our knowledge, the constructive proof we present below is not.
In fact, it is the key input for the degree-lowering isomorphisms of Section~\ref{sec:lowering-degree-general}, and  the blueprint for the $Q^n$-Algorithm of Section~\ref{sec: the-Q^n-Algorithm}.

For a topological Abelian group $G$, let $\mathcal F_G$ denote the \textit{sheaf of continuous $G$-valued maps} on $B$: $\mathcal F_G(V):=\mathrm{maps}(V,G)$ for open $V\subseteq B$, with the group structure inherited pointwise from $G$ and restriction maps given by restriction of functions. When $G$ is contractible, we have the following theorem.

\begin{theorem}[\Cech cohomology vanishes over $\mathcal{F}_E$]\label{Cech cohomology vanishes over soft sheaves}
Let $B$ be a paracompact Hausdorff space and $E$ a contractible topological Abelian group. Then $\check{H}^n(B;\mathcal{F}_E)=0$ for all $n>0$.
\end{theorem}

To prove Theorem~\ref{Cech cohomology vanishes over soft sheaves}, the key insight is that $\mathcal{F}_E$, together with the associated \Cech sheaves $\check{C}^n(\mathcal{U}, \mathcal{F}_E)$, are all \textit{soft sheaves} because of the contractibility of $E$; and, as we shall see, \Cech cohomology with coefficients in a soft sheaf vanishes in every positive degree, provided the underlying space is paracompact. Verifying that $\mathcal{F}_E$ and its associated \Cech sheaves are soft (Lemma~\ref{F_E is soft}), and then establishing this vanishing property for soft sheaves in general (Proposition~\ref{Cech Cohomology vanishes over F}), will together yield Theorem~\ref{Cech cohomology vanishes over soft sheaves}.

    \begin{definition}[Sections over arbitrary subspaces]
        Let $B$ be a topological space, and let $\mathcal{F}$ be a sheaf of Abelian groups over $B$. 
        For any subspace $X\subseteq B$, define the direct limit (colimit)
        $$\mathcal{F}(X):=\varinjlim_{X\subseteq U \open B} \mathcal{F}(U)$$
    \end{definition}

    \begin{definition}[Soft sheaves]
        Let $B$ be a paracompact space. 
        A sheaf $\mathcal{F}$ on $B$ is said to be \textit{soft} if $\forall K\closed B$, the restriction map $\mathcal{F}(B)\to \mathcal{F}(K)$ is surjective, i.e., every section of $\mathcal{F}$ over $K$ can be extended to a global section.
    \end{definition}

\begin{lemma}[Softness of $\mathcal{F}_E$ and its \Cech sheaves]
\label{F_E is soft}\label{Cech Sheaves are soft over F_E}
Let $B$ be a paracompact Hausdorff space, $E$ a contractible topological Abelian group, and $\mathcal{U}$ an open cover of $B$. Then the sheaf represented by $E$ over $B$, i.e., $\mathcal{F}_E:=\text{maps}(-,E)$, and the \Cech sheaves $\check{C}^n(\mathcal{U},\mathcal{F}_E)$ for all $n\geq 0$ are soft. 
\end{lemma}

\begin{proof}
The softness of $\mathcal F_E$ is due to Huber~\cite[Lemma~3.2]{huber1961homotopical}; we include a complete proof for the reader's convenience, and then extend the argument to  \Cech sheaves.

Let $h: E\times [0,1]\to E$ be a contracting homotopy, i.e., $h(-,1)=\mathbbm{1}_E$ and $h(-,0)\equiv p$ for some $p\in E$. Let $K\closed B$. We first treat $\mathcal{F}_E$. Let $[s]\in\mathcal{F}_E(K)$ be represented by $s:V\to E$ for some $K\subseteq V\open B$. Since $B$ is paracompact Hausdorff, hence normal, we may choose an open set $K\subseteq W\open B$ with $\overline{W}\subseteq V$, and a continuous map $g: B\to [0,1]$ with $g|_{\overline{W}}=1$ and $\operatorname{supp}(g)\subseteq V$, where $\operatorname{supp}(g):=\overline{\{b\in B\mid g(b)\neq 0\}}$ (Urysohn's lemma). Define $\tilde{s}: B\to E$ by
$$\tilde{s}(b) := \left\{ \begin{array}{ll}
h(s(b),g(b)) & b\in V,\\
p & b\notin V.
\end{array}\right.$$
Then $\tilde{s}$ is continuous, since it equals $h(s,g)$ on the open set $V$ and is constantly $p$ on the open set $B\setminus\operatorname{supp}(g)$, and these two open sets cover $B$ and the two descriptions agree on their overlap. Moreover $\tilde{s}|_W=s|_W$, so $\tilde{s}$ extends $[s]$.

For the \Cech sheaves, let $[s]\in\check{C}^n(\mathcal{U},\mathcal{F}_E)(K)$ be represented by
$$s=(s_{i_0\cdots i_n})_{i_0<\cdots<i_n}\in\prod_{i_0<\cdots<i_n}\text{maps}(U_{i_0\cdots i_n}\cap V, E)$$
for some $K\subseteq V\open B$, and choose $W$ and $g$ as above. For each $n$-simplex $(i_0,\dots,i_n)$, define $\tilde{s}_{i_0\cdots i_n}: U_{i_0\cdots i_n}\to E$ by
$$\tilde{s}_{i_0\cdots i_n}(b):= \left\{\begin{array}{ll}
    h(s_{i_0\cdots i_n}(b),g(b)) & b\in U_{i_0\cdots i_n}\cap V, \\
    p & b\in U_{i_0\cdots i_n}\setminus V.
\end{array}\right.$$
Each $\tilde{s}_{i_0\cdots i_n}$ is continuous, since it equals $h(s_{i_0\cdots i_n},g)$ on the open set $U_{i_0\cdots i_n}\cap V$ and is constantly $p$ on the open set $U_{i_0\cdots i_n}\setminus\operatorname{supp}(g)$, and these two open sets cover $U_{i_0\cdots i_n}$ and the two descriptions agree on their overlap. Since $g|_W\equiv 1$, we have $\tilde{s}_{i_0\cdots i_n}(b)=h(s_{i_0\cdots i_n}(b),1)=s_{i_0\cdots i_n}(b)$ for all $b\in W$, so $\tilde{s}:=(\tilde{s}_{i_0\cdots i_n})_{i_0<\cdots<i_n}$ extends $[s]$.
\end{proof}

We next establish two classical facts about soft sheaves on a paracompact space. They are the building blocks of Theorem~\ref{Cech cohomology vanishes over soft sheaves}, and thereby furnish the computational blueprint of the $Q^n$-Algorithm of Section~\ref{sec: the-Q^n-Algorithm}.

\begin{proposition}[Soft sheaf preserves exactness over global sections; cf.\ {\cite[II, Theorem~9.9]{bredonsheaf}}]
\label{Soft sheaf preserves exactness over global sections}
Let $B$ be a paracompact space and $0\to\mathcal A'\to\mathcal A\xrightarrow{\varphi}\mathcal A''\to0$ a short exact sequence of sheaves of Abelian groups on $B$ with $\mathcal A'$ soft. Then $0\to\mathcal A'(B)\to\mathcal A(B)\xrightarrow{\varphi_B}\mathcal A''(B)\to0$ is exact.
\end{proposition}

\begin{proof}
A complete proof can be found in \cite[II, Theorem~9.9]{bredonsheaf}; we outline the argument. Only the surjectivity of $\varphi_B$ is at issue. Given $s\in\mathcal A''(B)$, exactness of the sequence of sheaves provides local lifts $s_i\in\mathcal A(W_i)$ of $s|_{W_i}$ over a locally finite closed cover $\{W_i\}_{i\in I}$ of $B$. Suppose a lift $t$ of $s$ has been constructed over the closed set $F_J:=\bigcup_{j\in J}W_j$, and let $i\notin J$. The difference $s_i-t$ on $W_i\cap F_J$ maps to zero under $\varphi$, hence is a section of $\mathcal A'$; since $\mathcal A'$ is soft, it extends to some $\tilde s\in\mathcal A'(F_J)$, and then $t+\tilde s$ agrees with $s_i$ on $W_i\cap F_J$, so the two glue to a lift over $F_{J\cup\{i\}}$. A Zorn's lemma argument shows that this process reaches $F_I=B$; for a finite cover, finite induction suffices. This \emph{difference--extend--glue} step is precisely the inductive step of the $Q^n$-Algorithm (cf.\ Steps~1, 3, and~4 in the proof of Theorem~\ref{thm:Qn-correctness}).
\end{proof}

\begin{proposition}[Softness passes over on short exact sequences; cf.\ {\cite[II, Proposition~9.10]{bredonsheaf}}]
\label{Softness passes over on short exact sequences}
Let $B$ be a paracompact space and $0\to\mathcal A'\to\mathcal A\xrightarrow{\varphi}\mathcal A''\to0$ a short exact sequence of sheaves of Abelian groups over $B$ with $\mathcal A'$ soft. 
Then $\mathcal A$ is soft if and only if $\mathcal A''$ is soft.
\end{proposition}

\begin{proof}
For a closed set $K\subseteq B$, consider the commutative diagram of restriction maps
\begin{center}
\begin{tikzcd}
0 \arrow[r] & \mathcal A'(B) \arrow[r] \arrow[d,"(1)"] & \mathcal A(B) \arrow[r,"\varphi_B"] \arrow[d,"(2)"] & \mathcal A''(B) \arrow[r] \arrow[d,"(3)"] & 0\\
0 \arrow[r] & \mathcal A'(K) \arrow[r] & \mathcal A(K) \arrow[r,"\varphi_K"] & \mathcal A''(K) \arrow[r] & 0
\end{tikzcd}
\end{center}
The top row is exact by Proposition~\ref{Soft sheaf preserves exactness over global sections}, and the bottom row by the same argument applied over $K$, which is paracompact as a closed subspace. Since $\mathcal A'$ is soft, map~(1) is surjective, and a diagram chase shows that map~(2) is surjective if and only if map~(3) is.
\end{proof}

\begin{proposition}[\Cech cohomology vanishes over soft sheaves]
\label{Cech Cohomology vanishes over F}
Let $B$ be a paracompact space, $\mathcal F$ a soft sheaf on $B$, and $\mathcal U$ an open cover of $B$ such that $\check C^n(\mathcal U,\mathcal F)$ is soft for all $n\ge0$. Then $\check H^n(\mathcal U,\mathcal F)=0$ for all $n\ge1$.
\end{proposition}

\begin{proof}
Split the \Cech resolution of Lemma~\ref{Cech Resolution} into short exact sequences $0\to\ker d^n\to\check C^n(\mathcal U,\mathcal F)\to\ker d^{n+1}\to0$, where $\ker d^0:=\mathcal F$. Since $\mathcal F$ and all $\check C^n(\mathcal U,\mathcal F)$ are soft, Proposition~\ref{Softness passes over on short exact sequences} shows inductively that every $\ker d^n$ is soft, and Proposition~\ref{Soft sheaf preserves exactness over global sections} then shows that each of these sequences stays exact on global sections. Hence the global sections of the \Cech resolution form an exact sequence, i.e., $\check H^n(\mathcal U,\mathcal F)=0$ for $n\ge1$.
\end{proof}

\section{Methodology} \label{sec: methodology}

In this section, we present the methodology of the paper: how to compute, step by step, the Eilenberg-MacLane coordinate corresponding to a persistent cohomology class. Figure~\ref{fig:overview} gives an overview of this pipeline.

\begin{figure}[htbp]
\begin{tikzcd}
 PH^n(\mathcal{R}(X); \mathbb{Z}_q) \arrow[r, dashed, "\text{Select }\alpha"'] & H^n(\mathcal{R}_{2\alpha}(X);\mathbb{Z}_q) \arrow[r, "\ref{sec:bockstein-lifting}", "\text{Bockstein}"'] & H^n(\mathcal{R}_{2\alpha}(X);G) \arrow[r, "\ref{sec:from-simplicial-to-cech}","\text{Step 1}"'] & H^n(\check{\mathcal{C}}_\alpha(X);G) \arrow[d, "\ref{sec:from-simplicial-to-cech}", "\text{Step 2 \& 3}"'] \\
\check{H}^1(X^{(\alpha)};\mathcal{F}_{B^{n-1}G}) \arrow[d, "\ref{sec:chains-end}"] & \cdots \arrow[l, "\ref{sec:lowering-degree-general}","\cong"'] & \check{H}^{n-1}(X^{(\alpha)};\mathcal{F}_{BG}) \arrow[l, "\ref{sec:lowering-degree-general}","\cong"'] & \check{H}^n(X^{(\alpha)};\mathcal{F}_G) \arrow[l, "\ref{sec:lowering-degree-general}","\cong"'] \\
\check{H}^0(X^{(\alpha)};\mathcal{F}_{B^nG}) \arrow[r, "\ref{sec:chains-end}"] & {}[X^{(\alpha)}, B^nG] \arrow[r, equal] & {}[ X^{(\alpha)}, K(G,n)] &  \\
\end{tikzcd}
\caption{Overview of the pipeline from a persistent cohomology class to its corresponding Eilenberg-MacLane coordinate. Each arrow is labeled with the section carrying out that step.}
\label{fig:overview}
\end{figure}

We start with a persistent cohomology class $[\boldsymbol{\eta}]\in PH^n(\mathcal R(X);\mathbb Z_q)$, for a prime $q$, with persistence interval $[a,b)$ satisfying $2a<b$, and fix a scale $\alpha\in[a,\frac b2)$. The remainder of this section shows that the information carried by $[\eta_{2\alpha}]\in H^n(\mathcal R_{2\alpha}(X);\mathbb Z_q)$ can always be transported into a corresponding Eilenberg-MacLane coordinate $f_{\boldsymbol{\eta}}:X^{(\alpha)}\to K(G,n)$. Different choices of the coefficient group $G$ yield representations of different quality; this choice is addressed first, in Section~\ref{sec:bockstein-lifting}.

\subsection{Bockstein Lifting}\label{sec:bockstein-lifting}

Given a class $[\eta_{2\alpha}]\in H^n(\mathcal R_{2\alpha}(X);\mathbb Z_q)$, does it faithfully represent the underlying topological feature? Working with $\mathbb{Z}_q$ coefficients is standard practice for computational reasons — field coefficients guarantee that the persistence module decomposes into interval summands via the classification theorem (\ref{classification of pid}) — but the class $[\eta_{2\alpha}]$ obtained this way need not be intrinsically tied to $\mathbb{Z}_q$: the very same class could equally arise as the mod-$q$ reduction of a class with coefficients in $\mathbb Z$, in $\mathbb Z_{q^2}$, or more generally in any abelian group $G$ admitting a surjection onto $\mathbb{Z}_q$. More precisely, writing $K:=\mathcal R_{2\alpha}(X)$, the Bockstein long exact sequence for $0 \to \mathbb{Z} \to \mathbb{Z} \to \mathbb{Z}_q \to 0$ shows that a $\mathbb Z_q$-class is sensitive only to the free part of $H^n(K;\mathbb Z)$ and the $q$-primary torsion of $H^n(K;\mathbb Z)$ and $H^{n+1}(K;\mathbb Z)$. 
One therefore wishes to identify, among the coefficient groups $\mathbb Z$ and $\mathbb Z_{q^k}$, the one yielding the most faithful representation of the underlying feature.

The standard mechanism for recovering such a $G$, introduced by de Silva and Vejdemo-Johansson \cite{desilva2009persistent} and subsequently adopted throughout the topological-coordinates literature \cite{desilva2011persistent, perea2019sparse, schonsheck2024spherical, nordwald2026density}, proceeds by examining the Bockstein sequence associated to a short exact sequence of coefficients — a procedure we refer to as \textit{Bockstein lifting} throughout this paper. Given a short exact sequence of abelian groups $0\to H\xrightarrow{\iota} G\xrightarrow{\kappa} \mathbb{Z}_q\to 0$, the induced long exact sequence
$$\cdots \to H^n(K;H)\xrightarrow{\iota_*} H^n(K;G)\xrightarrow{\kappa_*} H^n(K;\mathbb{Z}_q)\xrightarrow{\Delta} H^{n+1}(K;H)\to \cdots$$
identifies the obstruction to lifting $[\eta_{2\alpha}]$ along $\kappa_*$ with its image $\Delta([\eta_{2\alpha}])\in H^{n+1}(K;H)$ under the connecting (Bockstein) homomorphism $\Delta$. When $\Delta([\eta_{2\alpha}])=0$, exactness guarantees a class $[\eta_G]\in H^n(K;G)$ with $\kappa_*([\eta_G])=[\eta_{2\alpha}]$, obtained concretely by correcting an integer representative of a cocycle for $\eta_{2\alpha}$ by a cochain divisible by $q$ (see Step~(4) of Section~\ref{sec: sparse-eilenberg-maclane-algorithm}). Such a lift is in general not unique: any two lifts differ by an element of $\iota_*(H^n(K;H))$, and the choice affects the quality of the resulting coordinate~\cite{gaukstad2025lifting}.

\begin{table}[htbp]
\centering
$$
\begin{array}{|c|c|}
\hline
k & 0\to H\xrightarrow{\iota}G\xrightarrow{\kappa}\mathbb Z_q\to 0\\
\hline
0 & 0\to 0\to\mathbb Z_q\to\mathbb Z_q\to 0\\
1 & 0\to\mathbb Z_q\xrightarrow{\cdot q}\mathbb Z_{q^2}\to\mathbb Z_q\to 0\\
2 & 0\to\mathbb Z_{q^2}\xrightarrow{\cdot q}\mathbb Z_{q^3}\to\mathbb Z_q\to 0\\
\cdots & 0 \to \cdots \to \cdots \to \cdots\to 0\\
\infty & 0\to\mathbb Z\xrightarrow{\cdot q}\mathbb Z\to\mathbb Z_q\to 0\\
\hline
\end{array}
$$
\caption{The tower of coefficient sequences tested by Bockstein lifting.}
\label{tab:bockstein-tower}
\end{table}

Since a lift to $\mathbb Z$ or to $\mathbb Z_{q^{k+1}}$ reduces modulo $q^k$ to a lift to $\mathbb Z_{q^k}$, the rows of Table~\ref{tab:bockstein-tower} form a chain of increasingly refined coefficient groups $\mathbb Z_q\leftarrow\mathbb Z_{q^2}\leftarrow\cdots\leftarrow\mathbb Z$, along which the rows whose obstruction vanishes form an initial segment; row $0$ always succeeds, since $H=0$. We take $G$ from the last row that succeeds, as this yields the most refined lift available; in particular, $G=\mathbb Z$ whenever row $\infty$ succeeds.

\subsection{From the Vietoris--Rips Complex to the $\alpha$-Offset} \label{sec:from-simplicial-to-cech}

Section~\ref{sec:bockstein-lifting} produced a class $[\eta_{G}]\in H^n(\mathcal R_{2\alpha}(X);G)$ on the Vietoris--Rips complex, on which persistent cohomology is computed in practice (e.g.\ by Ripser). Since the coordinate is to be defined on the $\alpha$-offset $X^{(\alpha)}$, we now transport this class to the \Cech cohomology of $X^{(\alpha)}$, to which the isomorphisms of Sections~\ref{sec:lowering-degree-general} and~\ref{sec:chains-end} apply.

\subsubsection*{1. The Rips--\Cech Interleaving}

The Vietoris--Rips and \Cech complexes (defined in \eqref{def: Rips and Cech Complexes}) are related by the inclusions~\cite{deSilvaGhrist2007}\footnote{When $X\subseteq\mathbb R^d$ with the Euclidean metric, the factor $2$ can be improved~\cite{deSilvaGhrist2007}.}
$$\mathcal R_\alpha(X)\ \subseteq\ \check{\mathcal C}_\alpha(X)\ \subseteq\ \mathcal R_{2\alpha}(X).$$
Restriction along the second inclusion maps $[\eta_G]$ to a class $[\eta_{\check{\mathcal C}}]\in H^n(\check{\mathcal C}_\alpha(X);G)$, which is nonzero: the composite $H^n(\mathcal R_{2\alpha}(X);G)\to H^n(\check{\mathcal C}_\alpha(X);G)\to H^n(\mathcal R_\alpha(X);G)$ is the persistent cohomology restriction map, and since both $\alpha$ and $2\alpha$ lie in the interval $[a,b)$, the first map in the factorization — under which $[\eta_G]$ maps to $[\eta_{\check{\mathcal C}}]$ — must itself be nonzero. This is what prevents $[\eta_G]$ from degenerating into the trivial class as it is transported to $H^n(\check{\mathcal C}_\alpha(X);G)$.

\subsubsection*{2. From Simplicial to \Cech Cohomology}

Write $B:=X^{(\alpha)}$ and $\mathcal U^\alpha:=\{B_\alpha(x)\}_{x\in X}$, and let $\mathcal C_G$ denote the \emph{constant presheaf} with value $G$ on $B$: $\mathcal C_G(U):=G$ for every nonempty $U\subseteq_{\text{open}} B$, with all restriction maps equal to $\mathbbm 1_G$, and $\mathcal C_G(\varnothing):=0$. A \Cech $n$-cochain of $\mathcal U^\alpha$ with coefficients in $\mathcal C_G$ assigns one element of $G$ to each nonempty intersection $U_{i_0\cdots i_n}=B_\alpha(x_{i_0})\cap\cdots\cap B_\alpha(x_{i_n})$. Since the $n$-simplices of $\check{\mathcal C}_\alpha(X)$ are precisely the tuples $\{x_{i_0},\ldots,x_{i_n}\}$ for which this intersection is nonempty, and all restriction maps are identities, these cochains and their coboundaries coincide term by term with the simplicial cochains of $\check{\mathcal C}_\alpha(X)$ and their coboundaries, with no further hypothesis. This yields the identification
$$H^n(\check{\mathcal C}_\alpha(X);G)\ =\ \check H^n(\mathcal U^\alpha,\mathcal C_G),$$
under which the class $[\eta_{\check{\mathcal C}}]$ of Step~1 corresponds to a class $[\eta_{\mathcal C_G}]\in\check H^n(\mathcal U^\alpha,\mathcal C_G)$.

\subsubsection*{3. Passing to the Sheaf $\mathcal F_G$}

Since the long exact sequence used in Section~\ref{sec:lowering-degree-general} requires a short exact sequence of \emph{sheaves}, we pass from $\mathcal C_G$ to its sheafification $\mathcal F_G$, the sheaf of continuous (equivalently, since $G$ is discrete, locally constant) $G$-valued maps. The natural map $\mathcal C_G\to\mathcal F_G$, regarding an element of $G$ as a constant map, induces a map on \Cech cohomology
$$\check H^n(\mathcal U^\alpha,\mathcal C_G)\ \to\ \check H^n(\mathcal U^\alpha,\mathcal F_G),$$
under which $[\eta_{\mathcal C_G}]$ corresponds to a class $[\eta_{\mathcal F_G}]$. If, in addition, the \emph{good cover} hypothesis holds — every nonempty intersection of the balls $B_\alpha(x)$ is connected (automatically satisfied, for instance, when $X\subseteq\mathbb R^d$ with the Euclidean metric) — then $\mathcal C_G$ and $\mathcal F_G$ agree on every nonempty intersection of $\mathcal U^\alpha$, since a locally constant $G$-valued function on a connected set is constant. In this case the map above is an isomorphism, so $[\eta_{\mathcal F_G}]$ is guaranteed to be nonzero.

Composing Steps~1--3 yields the class $[\eta_{\mathcal F_G}]\in\check H^n(\mathcal U^\alpha,\mathcal F_G)$, whose image in $\check H^n(B;\mathcal F_G)$ is the input to the isomorphisms $\check H^n(B;\mathcal F_G)\cong\check H^{n-1}(B;\mathcal F_{BG})\cong\cdots$ of Section~\ref{sec:lowering-degree-general}.

\subsection{Lowering the Degree of Čech Cohomology}\label{sec:lowering-degree-general}

From Section~\ref{sec:from-simplicial-to-cech}, we obtained a representative in $\check{H}^n(B,\mathcal{F}_G)$, where $B=X^{(\alpha)}$; our remaining task is to relate this group to $[B,B^nG]$. We aim to attain an isomorphism $\check{H}^n(B,\mathcal{F}_G)\cong [B,B^nG]$, and the entire construction reduces to a single building block, which we may apply recursively:
$$\check{H}^n(B,\mathcal{F}_G)\cong \check{H}^{n-1}(B,\mathcal{F}_{BG}),\qquad \forall n\geq 2$$
The idea behind this isomorphism is the long exact sequence in Čech cohomology induced by a short exact sequence of sheaves $0\to\mathcal F\to\mathcal G\to\mathcal H\to 0$, with the connecting homomorphism $\Delta$ linking cohomology across consecutive degrees. Crucially, whenever the middle sheaf $\mathcal G$ is chosen so that its higher Čech cohomology vanishes, the long exact sequence collapses into a family of isomorphisms $\Delta: \check{H}^n(B,\mathcal{H})\xrightarrow{\cong} \check{H}^{n+1}(B,\mathcal{F})$. Applying this observation to the short exact sequence induced by the universal bundle $G\to EG\to BG$, and iterating over the tower $G, BG, B^2G, ..., B^{n-1}G$ is precisely what will allow us to lower the degree of $\check{H}^n(B,\mathcal{F}_G)$ down to $\check{H}^1(B,\mathcal{F}_{B^{n-1}G})$, one bar construction at a time.

\begin{theorem}[Long exact sequence in \Cech cohomology; {\cite[IX.3.18]{miranda}}, {\cite[I.3.25]{SerreFAC}}]
\label{MirandaThm}
    Let $B$ be a paracompact Hausdorff space, and $0\to \mathcal{F} \xrightarrow{\alpha} \mathcal{G} \xrightarrow{\beta} \mathcal{H} \to 0$ a short exact sequence of sheaves of Abelian groups over $B$, where $\alpha$ and $\beta$ are maps of sheaves. 
    Then, there exist connecting homomorphisms $\Delta: \check{H}^n(B,\mathcal{H})\to \check{H}^{n+1}(B,\mathcal{F})$ for every $n\geq 0$ such that the sequence of \Cech cohomology groups
    $$0\to \check{H}^0(B,\mathcal{F})\to ...\to \check{H}^n(B,\mathcal{F})\xrightarrow{\alpha^*} \check{H}^n(B,\mathcal{G})\xrightarrow{\beta^*} \check{H}^n(B,\mathcal{H})\xrightarrow{\Delta} \check{H}^{n+1}(B,\mathcal{F})\to ...$$
is exact.
\end{theorem}

For Abelian topological groups $A, A'$, a continuous homomorphism $\phi: A\to A'$ induces, by postcomposition, a map of sheaves $\phi_\#:\mathcal F_A\to\mathcal F_{A'}$, $\phi_\#(f):=\phi\circ f$. In particular, the universal bundle $G\xrightarrow{i} EG\xrightarrow{j} BG$ of Section~\ref{section: a good model for K(G,n)} induces a sequence of sheaves
$$0\to \mathcal{F}_G\xrightarrow{i_\#} \mathcal{F}_{EG}\xrightarrow{j_\#} \mathcal{F}_{BG}\to 0.$$

\begin{lemma} \label{SES of sheaves} Let $G\xrightarrow{i}  EG\xrightarrow{j} BG$ be a universal principal $G$-bundle. 
Then, 
\[
0\to \mathcal{F}_G\xrightarrow{i_\#} \mathcal{F}_{EG}\xrightarrow{j_\#} \mathcal{F}_{BG}\to 0
\]
defined above is a short exact sequence of sheaves.
\end{lemma}

\begin{proof}
Surjectivity: given $s\in\mathcal F_{BG}(U)$ and $x\in U$, choose a local trivialization $\varphi: j^{-1}(W)\xrightarrow{\cong} W\times G$ with $s(x)\in W$; then $\tilde s:=\varphi^{-1}(s(\cdot),e)$ on $V:=s^{-1}(W)$ satisfies $j\circ\tilde s=s|_V$. Exactness in the middle and injectivity: $j\circ f\equiv e_{BG}$ means $f$ takes values in the fiber $j^{-1}(e_{BG})=i(G)$, and $i:G\to j^{-1}(e_{BG})$ is a homeomorphism (again by local triviality), so $\ker j_\#=i_\#(\mathcal F_G)\cong\mathcal F_G$.
\end{proof}

Now fix $n\geq 1$ and a discrete Abelian group $G$. Iterating the bar construction yields the tower of topological Abelian groups $G, BG, \dots, B^{n-1}G$ (Section~\ref{section: a good model for K(G,n)}), with universal bundles $B^kG\to EB^kG\to B^{k+1}G$ for $0\le k\le n-2$, where $B^0G:=G$. By Lemma~\ref{SES of sheaves} and Theorem~\ref{MirandaThm}, each of them yields a long exact sequence
$$\cdots \to \check{H}^i(B, \mathcal{F}_{EB^kG})\to \check{H}^i(B, \mathcal{F}_{B^{k+1}G})\xrightarrow{\Delta} \check{H}^{i+1}(B, \mathcal{F}_{B^kG})\to \check{H}^{i+1}(B, \mathcal{F}_{EB^kG})\to \cdots.$$
Since $EB^kG$ is contractible (Proposition~\ref{Contractibility of EG}), Theorem~\ref{Cech cohomology vanishes over soft sheaves} gives $\check{H}^i(B, \mathcal{F}_{EB^kG})=0$ for all $i>0$, so the connecting homomorphism $\Delta$ is an isomorphism for every $i\geq 1$. For $1\le i\le n-1$, let $\Delta_i$ denote this isomorphism in degree $i$ for $k=n-1-i$:
$$\Delta_i: \check{H}^i(B,\mathcal{F}_{B^{n-i}G}) \xrightarrow{\cong} \check{H}^{i+1}(B,\mathcal{F}_{B^{n-1-i}G}).$$
These compose to the chain
$$\check{H}^1(B,\mathcal{F}_{B^{n-1}G}) \xrightarrow[\cong]{\Delta_{1}} \check{H}^2(B,\mathcal{F}_{B^{n-2}G}) \xrightarrow[\cong]{\Delta_{2}}\cdots  \xrightarrow[\cong]{\Delta_{n-2}} \check{H}^{n-1}(B,\mathcal{F}_{BG}) \xrightarrow[\cong]{\Delta_{n-1}} \check{H}^n(B,\mathcal{F}_{G}).$$
Read from right to left, this chain lowers $\check H^n(B;\mathcal F_G)$ to $\check H^1(B;\mathcal F_{B^{n-1}G})$; Section~\ref{sec:chains-end} treats the remaining step from degree one to homotopy classes of maps $B\to B^nG$.

\subsection{The Chain's End: From $\check{H}^0$ to Homotopy Classes} \label{sec:chains-end}

The isomorphism chain established in Section~\ref{sec:lowering-degree-general} reduces the problem of computing 
$\check{H}^n(B,\mathcal{F}_G)$ to that of computing $\check{H}^1(B,\mathcal{F}_H)$, where $H:=B^{n-1}G$; it remains to relate the latter to   homotopy classes of maps $[B,BH]=[B,B^nG]$. 
Notably, the vanishing theorem underlying~\ref{sec:lowering-degree-general} (Theorem~\ref{Cech cohomology vanishes over soft sheaves}) only forces $\check{H}^i(B,\mathcal{F}_{EH})=0$ for $i>0$; at $i=0$, the group $\check{H}^0(B,\mathcal{F}_{EH})$ need not vanish, so the degree-lowering mechanism differs slightly at the chain's end, and a separate argument is needed to close the chain.

The long exact sequence
$$\cdots \to \check{H}^0(B,\mathcal{F}_{EH}) \xrightarrow{j_\#^*} \check{H}^0(B,\mathcal{F}_{BH}) \xrightarrow{\Delta} \check{H}^1(B,\mathcal{F}_{H}) \xrightarrow{i_\#^*} \check{H}^1(B,\mathcal{F}_{EH}) \to \cdots$$
along with  $\check{H}^1(B,\mathcal{F}_{EH})=0$, implies that this long exact sequence yields the surjectivity of $\Delta$ (but NOT bijectivity), hence the preimage of any representative in $\check{H}^1(B,\mathcal{F}_H)$ still exists. At the same time, $\check{H}^0(B,\mathcal{F}_{BH})\cong \mathcal{F}_{BH}(B)$ through the gluing axiom of sheaves, and $\mathcal{F}_{BH}(B)$ — the group of continuous maps $B\to BH$ under pointwise addition — surjects onto $[B,BH]:=\mathcal{F}_{BH}(B)/\simeq$ recovered as the quotient of $\mathcal{F}_{BH}(B)$ under the homotopy relation. To identify $\check{H}^1(B,\mathcal{F}_H)$ with $[B,BH]$, it therefore suffices to show that $\Delta$ and the quotient map onto $[B,BH]$ share the same fibers — that is, that $\Delta(f_0)=\Delta(f_1)$ if and only if $f_0\simeq f_1$ for any $f_0,f_1\in \mathcal{F}_{BH}(B)$. This is where the argument genuinely requires the geometry of 
$EH$, rather than the exactness of the sequence alone.

\begin{center}
\begin{tikzcd}
{\check{H}^1(B,\mathcal{F}_H)} \arrow[r, "?", dashed] & {[B,BH]}\\
{\check{H}^0(B,\mathcal{F}_{BH})} \arrow[u, "\Delta", two heads] \arrow[r, "\cong", shift left] & \mathcal{F}_{BH}(B) \arrow[l] \arrow[u, "/\simeq", two heads]
\end{tikzcd}
\end{center}

\begin{lemma}[Homotopy via Liftability]
\label{homotopy iff liftable}
For any $f_0,f_1\in \mathcal{F}_{BH}(B)$, $\Delta(f_0)=\Delta(f_1)$ if and only if $f_0\simeq f_1$.
\end{lemma}

\begin{proof} $(\Rightarrow)$
Since $\Delta$ is a homomorphism, $\Delta(f_0)=\Delta(f_1)$ iff $\Delta(f_1-f_0)=0$, iff $f_1-f_0\in \ker \Delta=\text{im} j_\#^*$ by exactness. Therefore, there exists continuous map $g: B\to EH$ such that $f_1-f_0=j_\#^*(g)=j\circ g$.

Since $EH$ is contractible (Proposition~\ref{Contractibility of EG}), $g$ is null-homotopic via $h: EH\times [0,1]\to EG$ with $h(-,0)=e_{EH}, h(-,1)=\mathbbm{1}_{EH}$. Construct $\overline{f}: B\times [0,1]\to BH$ by $\overline{f}(b,t) := j(h(g(b),t))+f_0(b)$. Then, $\overline{f}$ is continuous since all components defining it is continuous; Meanwhile, $\overline{f}(b,0)=j(h(g(b),0))+f_0(b)=j(e_{EH})+f_0(b)=e_{BH}+f_0(b)=f_0(b)$, and $\overline{f}(b,1)=j(h(g(b),1))+f_0(b)=j(g(b))+f_0(b)=(f_1-f_0)(b)+f_0(b)=f_1(b)$, so $\overline{f}$ gives an explicit homotopy between $f_0$ and $f_1$.

$(\Leftarrow)$ Consider $f_0\simeq f_1$, then there exists a continuous map $\overline{f}:B\times[0,1]\to BH$ with $\overline{f}(-,0)=f_0$ and $\overline{f}(-,1)=f_1$. Define $\varphi: B\times[0,1]\to BH$ by $\varphi(b,t):=\overline{f}(b,t)-f_0(b)$, which is continuous since $\overline{f}$ and $f_0$ are, and satisfies $\varphi(b,0)=f_0(b)-f_0(b)=e_{BH}$ and $\varphi(b,1)=f_1(b)-f_0(b)=(f_1-f_0)(b)$.

Since $B$ is paracompact Hausdorff, the principal bundle $H\to EH\xrightarrow{j} BH$ is numerable, hence a Hurewicz fibration. The constant map $\Psi_0:B\to EH$, $\Psi_0(b):=e_{EH}$, satisfies $j(\Psi_0(b))=j(e_{EH})=e_{BH}=\varphi(b,0)$, so $\Psi_0$ is a lift of $\varphi(-,0)$ through $j$. By the homotopy lifting property applied to $\varphi$ with initial lift $\Psi_0$, there exists a continuous map $\Phi:B\times[0,1]\to EH$ such that $\Phi(-,0)=\Psi_0$ and $j\circ\Phi=\varphi$.

Set $g:B\to EH$ by $g(b):=\Phi(b,1)$. Then $j(g(b))=j(\Phi(b,1))=\varphi(b,1)=(f_1-f_0)(b)$, so $j\circ g=f_1-f_0$, i.e., $f_1-f_0=j_\#^*(g)\in\text{im } j_\#^*=\ker\Delta$ by exactness. Hence $\Delta(f_1-f_0)=0$, and since $\Delta$ is a homomorphism, $\Delta(f_0)=\Delta(f_1)$.
\end{proof}

Lemma~\ref{homotopy iff liftable} shows that $\Delta$ and the quotient map onto $[B,BH]$ share the same fibers, so $\Delta$ descends to a bijection $$\overline{\Delta}: [B,BH] \xrightarrow{\cong} \check{H}^1(B,\mathcal{F}_H)$$ which is in fact a group isomorphism, since $\Delta$ is a homomorphism and the group structure on $[B,BH]$ is induced by the (continuous) addition on $BH$. Combining this with $H=B^{n-1}G$ and the isomorphism chain of Section 3.4 yields the master isomorphism:

$$\check{H}^n(B,\mathcal{F}_G)\cong \check{H}^{n-1}(B,\mathcal{F}_{BG}) \cong \cdots \cong \check{H}^1(B,\mathcal{F}_{B^{n-1}G}) \xrightarrow[\cong]{\overline{\Delta}} [B,B^nG]=[B,K(G,n)]$$

\subsection{Milnor variant: Relation to the Classical Principal-Bundle Approach}\label{sec: Milnor-variant}

Sections~\ref{sec:bockstein-lifting}--\ref{sec:lowering-degree-general} together furnish a complete pipeline: starting from a persistent cohomology class $[\boldsymbol{\eta}]$ satisfying the appropriate conditions, together with a choice of scale $\alpha$, one obtains a map $f_\eta: X^{(\alpha)}\to B^nG$. This is referred to as the \emph{(sparse) Eilenberg-MacLane coordinate} generated by $\eta$ and $\alpha$.

In this subsection we present a variant of this pipeline, called the \emph{Milnor Variant}, grafted onto the classical approach of de Silva, Morozov, and Vejdemo-Johansson~\cite{desilva2011persistent} and Perea~\cite{perea2019sparse}. This variant produces instead a map $X^{(\alpha)}\to\mathcal BB^{n-1}G$, where $\mathcal BH$ denotes the Milnor construction model of a classifying space (the details are clarified later in this section). The idea is to carry out Sections~\ref{sec:bockstein-lifting}--\ref{sec:lowering-degree-general} exactly as in the canonical pipeline, arriving at a class in $\check H^1(X^{(\alpha)},\mathcal F_{B^{n-1}G})$, but replaces the final step by the classical correspondence ($B:=X^{(\alpha)}$)
\begin{equation}
\check{H}^1(B, \mathcal{F}_{H})\ \cong\ \mathrm{Prin}_H(B)\ \cong\ [B, \mathcal{B}H]
\label{eq:milnor-correspondence}
\end{equation}

The chief advantage of this variant pipeline is twofold: it not only connects the canonical pipeline to the classical principal-bundle approach, but also admits better geometric properties. Unlike the iterated bar construction, the Milnor construction passes a metric from the coefficient group to its classifying space: whenever $B^{n-1}G$ carries a right-invariant metric, $\mathcal B(B^{n-1}G)$ inherits one. This in turn makes distance-dependent procedures --- such as Isomap, multidimensional scaling (MDS), or a further persistent cohomology computation --- directly applicable to the resulting coordinate, and thereby helps visualize the coordinate space. The metric must still originate in $B^{n-1}G$, which the bar construction does not supply in general; the advantage is realized whenever $B^{n-1}G$ can be replaced by a metric group. Examples include $\mathcal B(B\mathbb Z)$, via $B\mathbb Z \simeq S^1$; more generally $\mathcal B(B\mathbb Z^k)$, via $B\mathbb Z^k \simeq T^k$; and $\mathcal B(B\mathbb Z_q)$, when the coordinate stays within the 1-skeleton of $B\mathbb Z_q$, as in Experiment~5.

\subsubsection*{How does correspondence (\ref{eq:milnor-correspondence}) work} Recall the Milnor construction of a classifying space \cite[Def.~4.11.1]{husemollerfiberbundles}: for a topological group $H$, one sets $\mathcal{E}H:=\bigsqcup_{n\geq 0}H^{n+1}\times\Delta_n/\sim$
with the equivalence relation $\sim$ given by replacing an entry $h_i$ arbitrarily whenever $t_i=0$. Unlike the bar construction, $H$ acts \emph{diagonally} on every coordinate simultaneously: $[(h_0,\ldots,h_n),(t_0,\ldots,t_n)]\cdot h:=[(h_0h,\ldots,h_nh),(t_0,\ldots,t_n)],$
and $\mathcal{B}H:=\mathcal{E}H/H$. Composing the two correspondences in (\ref{eq:milnor-correspondence}) yields an explicit formula: given $\eta\in\check Z^1(B,\mathcal F_H)$ and a partition of unity $\{\phi_k: B\to [0,1]\}$ subordinated to an open cover $\mathcal{U}$ of $B$, the map $f_\eta: B\to \mathcal{B}H$ defined by
\begin{equation}\label{eq:milnor-formula}
f_\eta(b):=\big[(\eta_{j\sigma_0}(b),\dots,\eta_{j\sigma_l}(b)),(\phi_{\sigma_0}(b),\dots,\phi_{\sigma_l}(b))\big],
\end{equation}
where $\sigma_b:=\{k\mid b\in U_k\}=\{\sigma_0,\dots,\sigma_l\},\ j\in\sigma_b$, is a classifying map of the principal $H$-bundle determined by $\eta$ \cite[Thm.~4.12.2]{husemollerfiberbundles}, \cite{perea2019sparse}. It is independent of the choice of $j$: by the cocycle condition, replacing $j$ by $j'$ shifts every entry by $\eta_{j'j}(b)$, which is absorbed by the diagonal action. This is the formula implemented in Algorithm~\ref{alg: Milnor-Variant}.

\begin{remark}[A Lighter-Weight Construction]
Unlike the Milnor variant which requires a partition of unity $\{\phi_j\}$, the canonical pipeline uses bump functions $g_i$ subject only to local boundary conditions, with no constraint relating $g_i$ to $g_j$ for $i\neq j$. In general, such bump functions exist on any normal space by Urysohn's lemma, whereas partitions of unity require paracompactness.
\end{remark}

\section{Algorithm} \label{sec:algorithm}

This section carries out the precise computations underlying each arrow of the pipeline outlined at the beginning of Section \ref{sec: methodology}, assembling them into what we call the \emph{Eilenberg-MacLane Algorithm}. Our central task is to make explicit the degree-lowering step of Section~\ref{sec:lowering-degree-general}, which rests on Theorem~\ref{Cech cohomology vanishes over soft sheaves} ($\check H^i(B,\mathcal F_{EG})=0$ for $i>0$).

We start with the following general \emph{coboundary lifting} question: given $\eta\in\check C^n(\mathcal U,\mathcal F_{EG})(B)$ with $d^n_B(\eta)=0$, can one algorithmically compute a $\tau\in\check C^{n-1}(\mathcal U,\mathcal F_{EG})(B)$ such that $d^{n-1}_B(\tau)=\eta$? 
The answer lies within the proof of Theorem~\ref{Cech cohomology vanishes over soft sheaves} itself, and requires assembling the proofs into a single, explicit construction. 
As the core difficulty and main contribution of this paper, the work of computing such a $\tau$ is referred to, throughout this paper, as the \emph{$Q^n$-Algorithm}.

\begin{center}
\begin{tikzcd}
    \check{C}^{n-1}(\mathcal{U},\mathcal{F}_G)(B) \arrow[r,"i_{\#}"] \arrow[d,"d^{n-1}_B"] & \check{C}^{n-1}(\mathcal{U},\mathcal{F}_{EG})(B)\arrow[d,bend left,"d^{n-1}_B"] \arrow[r,"j_{\#}"] & \check{C}^{n-1}(\mathcal{U},\mathcal{F}_{BG})(B)\arrow[d,"d^{n-1}_B"] \\
    \check{C}^n(\mathcal{U},\mathcal{F}_G)(B) \arrow[r,"i_{\#}"] & \check{C}^n (\mathcal{U},\mathcal{F}_{EG})(B) \arrow[u,bend left,"Q^n"] \arrow[r,"j_{\#}"] & \check{C}^n(\mathcal{U},\mathcal{F}_{BG})(B) \\
\end{tikzcd}
\end{center}

\subsection{The $Q^n$-Algorithm: From $\ker d^n_B$ to $\check C^{n-1}(\mathcal U,\mathcal F_{EG})(B)$} \label{sec: the-Q^n-Algorithm}

Suppose we are given $\eta\in\ker d^n_B$. Throughout, $B:=X^{(\alpha)}$, and $\mathcal U:=\{U_0,U_1,\ldots,U_{N-1}\}$ is a finite open cover of $B$ consisting of $N$ open sets; in practice, the most natural choice is the \textit{union of balls}, $\mathcal U:=\{B_\alpha(x)\}_{x\in X}$.

We first construct an \textit{open cover filtration} ($N$ open covers in total)
\[\mathcal{U}=\mathcal{U}^{\gamma_0}\to \mathcal{U}^{\gamma_1}\to \cdots \to \mathcal{U}^{\gamma_{N-1}}\]
where each arrow denotes a shrinking of the preceding cover — such a shrinking is guaranteed to exist for any open cover of a paracompact Hausdorff space by the Shrinking Lemma~\cite[Lemma~41.6]{munkres2000topology}. If $\mathcal U$ is chosen as a union of balls, such a filtration may be constructed explicitly: for any $r>0$, set $\mathcal U^r:=\{B_r(x)\}_{x\in X}$; then any strictly decreasing sequence of radii $\alpha=\gamma_0>\gamma_1>\cdots>\gamma_{N-1}$, chosen so that $\mathcal U^{\gamma_{N-1}}$ still covers $B$, yields the desired open cover filtration $\mathcal U=\mathcal U^{\gamma_0}\to\mathcal U^{\gamma_1}\to\cdots\to\mathcal U^{\gamma_{N-1}}$.

Next, for $r>0$, write $\tilde U_i^r:=U_0^r\cup U_1^r\cup\cdots\cup U_i^r$ for brevity. We must construct $N-1$ bump functions
$$g_1,g_2,\ldots,g_{N-1}:B\to[0,1],$$
Each $g_i$ is required to be identically $1$ on $\overline{U_i^{\gamma_i}}\cap \overline{\tilde U_{i-1}^{\gamma_i}}$ (overline means the closure), and identically $0$ on $U_i^{\gamma_0}\setminus\bigl(U_i^{\gamma_0}\cap\tilde U_{i-1}^{\gamma_{i-1}}\bigr)$, from \ref{Cech Sheaves are soft over F_E}. For $\mathcal{U}$ chosen as a union of balls, an evident choice is to use piecewise-linear functions to accomplish this: define
$$\Phi_i^{b>a}:\mathbb R_{\geq 0}\to[0,1];\qquad \Phi_i^{b>a}(x):=\begin{cases}
1 & \text{if } 0\leq x\leq a\\[2pt]
\dfrac{b-x}{b-a} & \text{if } a< x\leq b\\[6pt]
0 & \text{if } x> b
\end{cases}$$
and set
$$g_i(x):=\min\Bigl\{\max_{j=0}^{i-1}\Phi_j^{\gamma_{i-1}>\gamma_i}\bigl(d(x,x_j)\bigr),\ \Phi_i^{\gamma_0>\gamma_i}\bigl(d(x,x_i)\bigr)\Bigr\}.$$
The reader may verify that $g_i$, so defined, satisfies the required conditions.
\begin{figure}[htbp]
    \centering
    \fbox{\includegraphics[width=0.55\textwidth]{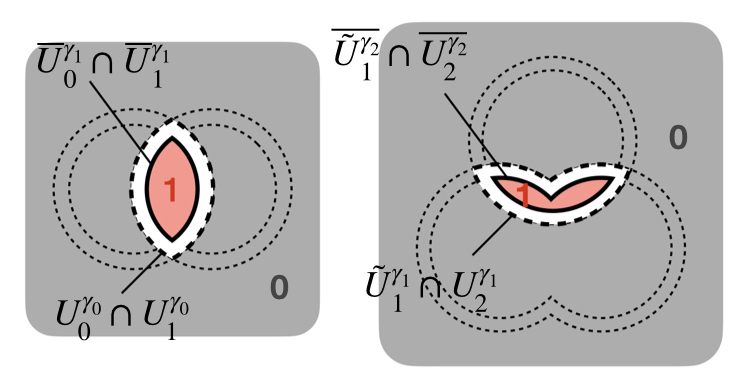}}
    \fbox{\includegraphics[width=0.25\textwidth]{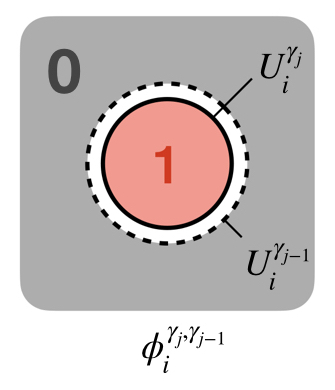}}
    \caption{Left: Regions in which the bump functions $g_1$ and $g_2$ take the value $0$ (gray) and $1$ (red), illustrated for the case where $\mathcal U$ consists of 3 open sets. Dashed boundaries denote open boundaries. Right: Regions in which $\Phi_i^{\gamma_{j-1}>\gamma_j}$ takes the value $0$ (gray) and $1$ (red).}
\end{figure}

Before turning to the explicit construction of $\tau$, it is worth explaining why bump functions of this particular shape are needed. The recursive construction of Section~4.1 proceeds by gluing a growing partial preimage, one open set at a time: having already assembled a section over $\tilde U_{i-1}$, the next step must incorporate the new information carried over $U_i$, while agreeing with what has already been glued on the overlap $U_i\cap\tilde U_{i-1}$. The function $g_i$ is precisely the interpolating device that mediates this transition: it must equal $1$ on the region where the new and previously-assembled data overlap (so that the two are blended together there), and it must equal $0$ on the region of $U_i$ not yet touched by $\tilde U_{i-1}$ (so that the newly introduced data is left unmodified there). The shrinking sequence $\gamma_0>\gamma_1>\cdots>\gamma_{N-1}$ supplies exactly the buffer needed to construct such a function continuously, and the piecewise-linear choice above is simply the most direct way to realize this behavior for a cover by balls.

Next, we will introduce three basic operations in the $Q^n$-algorithm.\\

\textit{1. Lifting}

We first define the \emph{lift} of a cocycle. For an $n$-\Cech cocycle $\theta\in(\ker d^n)(V)$, i.e., $\theta\in\check C^n(\mathcal U,\mathcal F_{EG})(V)$ such that $d_V^n(\theta)=0$, define its lift
$$\mathcal L(\theta):=\bigl(\mathcal L(\theta)^{(0)},\ldots,\mathcal L(\theta)^{(N-1)}\bigr)\in\prod_{i=0}^{N-1}\check C^{n-1}(\mathcal U,\mathcal F_{EG})(V\cap U_i)$$
by
$$\mathcal L(\theta)^{(i)}_{i_0\ldots i_{n-1}}:=\begin{cases}
\text{the constant map at } e_{EG}, & \text{if } i\in\{i_0,\ldots,i_{n-1}\}\\
\theta_{ii_0\ldots i_{n-1}}, & \text{if } i\notin\{i_0,\ldots,i_{n-1}\}.
\end{cases}$$

\textit{2. Extension}

We next define the \emph{extension} of a cochain over a single open set. Recall the argument used in the proof of Lemma~\ref{Cech Sheaves are soft over F_E}. Suppose we are given a cochain $\theta\in\check C^n(\mathcal U,\mathcal F_{EG})(V)$, a closed set $K\subseteq_{\text{closed}}B$ with $K\subseteq V$, an open set $W\supseteq V$, and a function $g:B\to[0,1]$ such that $g|_K\equiv1$ and $g|_{W-V}\equiv0$.
\footnote{Strictly speaking, continuity of $\tilde\theta$ at $\partial V$ for arbitrary $\theta$ requires $\operatorname{supp}(g)\cap W\subseteq V$, slightly stronger than $g|_{W\setminus V}\equiv 0$. For the bump functions $g_i$ above, this is achieved by shrinking $\gamma_{i-1}$ in the first factor by an arbitrarily small amount; we suppress this for readability.}
Then define $\tilde\theta\in\check C^n(\mathcal U,\mathcal F_{EG})(W)$ by
$$\tilde\theta_{i_0\ldots i_n}(x):=\begin{cases}
h\bigl(\theta_{i_0\ldots i_n}(x),g(x)\bigr), & \text{if } x\in V\cap U_{i_0\ldots i_n},\\
e_{EG}, & \text{if } x\in(W-V)\cap U_{i_0\ldots i_n},
\end{cases}$$
where $h$ is the contracting homotopy of $EG$ from Proposition~\ref{Contractibility of EG}. We call $\tilde\theta$ the \emph{extension of $\theta$ from $V$ to $W$ via $(g,K)$}, and write $\tilde\theta=\mathbf{Ext}^W_V\bigl(\theta;(g,K)\bigr)$. The underlying $K$ here indicates that it is guaranteed that $\tilde{\theta}|_K = \theta$.\\

\textit{3. Gluing}

Finally, we define the \emph{gluing} of two agreeing sections. Suppose $V_1,V_2\subseteq_{\text{open}}B$, and $\theta\in\check C^{n}(\mathcal U,\mathcal F_{EG})(V_1)$, $\theta'\in\check C^{n}(\mathcal U,\mathcal F_{EG})(V_2)$ satisfy $\theta|_{V_1\cap V_2}=\theta'|_{V_1\cap V_2}$. Define $\mathbf{Glue}(\theta,\theta')\in\check C^{n}(\mathcal U,\mathcal F_{EG})(V_1\cup V_2)$ by
$$\mathbf{Glue}(\theta,\theta')_{i_0\ldots i_{n}}(x):=\begin{cases}
\theta_{i_0\ldots i_{n}}(x), & \text{if } x\in V_1\cap U_{i_0\ldots i_{n}},\\
\theta'_{i_0\ldots i_{n}}(x), & \text{if } x\in (V_2\backslash V_1)\cap U_{i_0\ldots i_{n}}.
\end{cases}$$

\begin{algorithm}[H] \label{alg: Q^n-Algorithm}
\caption{$Q^n: \ker d^n_B \to \check C^{n-1}(\mathcal U,\mathcal F_{EG})(B)$}
	\begin{algorithmic}[1]
        \Require $\eta = (\eta_{i_0...i_n})_{i_0<...<i_n} \in \ker d_B^n\subseteq \check{C}^n(\mathcal{U}, \mathcal{F}_{EG})(B)$
        \Require Open cover filtration: $\mathcal U=\mathcal U^{\gamma_0}\to\mathcal U^{\gamma_1}\to\cdots\to\mathcal U^{\gamma_{N-1}}$ $(N=\# \mathcal{U})$
        \Require Bump functions $g_1, g_2,..., g_{N-1}: B\to [0,1]$
            \State $\eta^{(i)}:= \mathcal{L}(\eta)^{(i)}$, for $i=0,1,..., N-1$
            \State Set $t^{(0)}:= \eta^{(0)}$
		\For {$i=1,2,...,N-1$}
                \State $\tilde{t}^{(i)} := t^{(i-1)}|_{U_i^{\gamma_0}\cap \tilde{U}_{i-1}^{\gamma_{i-1}}}-\eta^{(i)}|_{U_i^{\gamma_0}\cap \tilde{U}_{i-1}^{\gamma_{i-1}}}$
			    \State $\bar{t}^{(i)} := d^{n-2}_{U_i^{\gamma_0}}\bigl(\mathbf{Ext}^{U_i^{\gamma_0}}_{U_i^{\gamma_0}\cap \tilde{U}_{i-1}^{\gamma_{i-1}}}\bigl(\mathcal{L}(\tilde{t}^{(i)})^{(i)}; (g_i, U_i^{\gamma_i}\cap \tilde{U}_{i-1}^{\gamma_i})\bigr)\bigr)$
                \State $t^{(i)}:= \mathbf{Glue}\bigl( t^{(i-1)}|_{\tilde{U}_{i-1}^{\gamma_i}}, (\eta^{(i)}+\bar{t}^{(i)})|_{U_i^{\gamma_i}}\bigr)$
		  \EndFor
        \State Set $Q^n(\eta) := t^{(N-1)} \in \check{C}^{n-1}(\mathcal{U}, \mathcal{F}_{EG})(B)$.
        \State As a consequence of Theorem~\ref{thm:Qn-correctness}, $d^{n-1}(Q^n(\eta))=\eta$.
	\end{algorithmic}
\end{algorithm}

\begin{figure}[htbp]
    \centering
    \fbox{\includegraphics[width=0.53\textwidth]{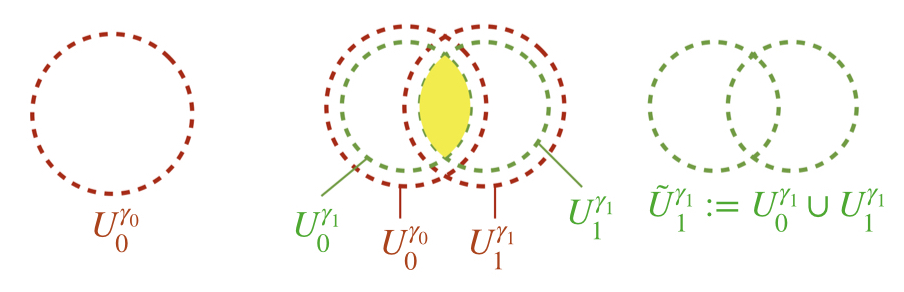}}
    \fbox{\includegraphics[width=0.41\textwidth]{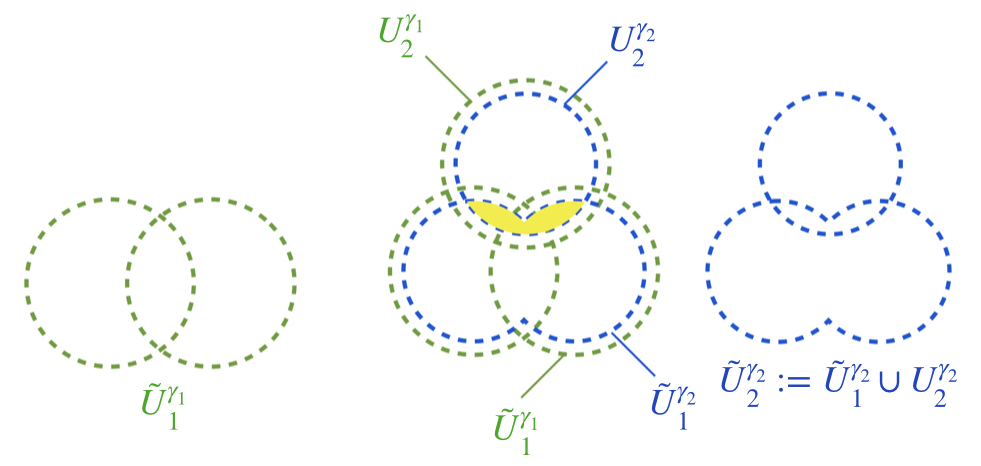}}
    \caption{The section of the \Cech cochain constructed by the $Q^n$-algorithm at $i=1$ and $i=2$. Through the operations of taking the difference, extending, and gluing, the domain grows step by step from $U_0^{\gamma_0}$ to $\tilde U_1^{\gamma_1}$, and then from $\tilde U_1^{\gamma_1}$ to $\tilde U_2^{\gamma_2}$.}
    \label{fig:qn-recursion-domains}
\end{figure}

The idea behind Algorithm~\ref{alg: Q^n-Algorithm} is as follows. A preimage of $\eta$ under $d^{n-1}$ always exists locally: on each $U_i$, the lift $\eta^{(i)}=\mathcal L(\eta)^{(i)}$ satisfies $d^{n-1}\eta^{(i)}=\eta$. The algorithm assembles these local preimages one open set at a time. Having built a preimage $t^{(i-1)}$ on $\tilde U_{i-1}^{\gamma_{i-1}}$, the new piece $\eta^{(i)}$ generally disagrees with it on their overlap. Their difference $\tilde t^{(i)}$ is a cocycle, hence locally a coboundary, and the contracting homotopy of $EG$ allows the resulting correction term $\bar t^{(i)}$ to be extended from the overlap to all of $U_i^{\gamma_0}$. After this correction, $\eta^{(i)}+\bar t^{(i)}$ agrees with $t^{(i-1)}$ on a slightly smaller overlap, where the two are glued; the decreasing radii $\gamma_0>\gamma_1>\cdots$ provide exactly this margin. Figure~\ref{fig:qn-recursion-domains} illustrates how the domain grows, and Theorem~\ref{thm:Qn-correctness} makes the argument precise.

\begin{theorem}[Correctness of the $Q^n$-Algorithm]\label{thm:Qn-correctness}
Let $B$ be a metric space, $\mathcal U=\{U_0,\dots,U_{N-1}\}$ a finite cover of $B$ by open balls $U_i=U_i^{\gamma_0}$, and $\gamma_0>\gamma_1>\cdots>\gamma_{N-1}$ radii such that $\mathcal U^{\gamma_{N-1}}$ still covers $B$. For $1\le i\le N-1$ put $V_i:=U_i^{\gamma_0}\cap\tilde U_{i-1}^{\gamma_{i-1}}, K_i:=\overline{U_i^{\gamma_i}}\cap\overline{\tilde U_{i-1}^{\gamma_i}},$ and suppose the bump functions $g_i: B\to [0,1]$ satisfy
\begin{enumerate}
\item[(B1)] $g_i\equiv 1$ on $K_i$;
\item[(B2)] $\operatorname{supp}(g_i)\cap U_i^{\gamma_0}\subseteq V_i$ (cf.\ the footnote in the definition of $\mathbf{Ext}$).
\end{enumerate}
Then for every $\eta\in\ker d^n_B$ and every $0\le i\le N-1$,
$$t^{(i)}\in\check C^{n-1}(\mathcal U,\mathcal F_{EG})\big(\tilde U_i^{\gamma_i}\big)\quad\text{and}\quad d^{n-1}t^{(i)}=\eta|_{\tilde U_i^{\gamma_i}}.$$
In particular $Q^n(\eta)\in\check C^{n-1}(\mathcal U,\mathcal F_{EG})(B)$ and $d^{n-1}\big(Q^n(\eta)\big)=\eta$. (For $n=1$ we use the convention $\check C^{-1}(\mathcal U,\mathcal F_{EG}):=\mathcal F_{EG}$ and $d^{-1}:=e$ used in Lemma~\ref{Cech Resolution}.)
\end{theorem}

\begin{proof}
We use two facts. First, for any open $W\subseteq U_i$ and any $\theta\in(\ker d^m)(W)$, the computation in the proof of Lemma~\ref{Local fineness of Cech Cohomology} applies verbatim with $B$ replaced by $W$, and gives
\begin{equation*}
d^{m-1}\,\mathcal L(\theta)^{(i)}=\theta. \tag{$\ast$}
\end{equation*}
Second, the \Cech coboundary on $\check C^\bullet(\mathcal U,\mathcal F_{EG})$ is computed pointwise, so two cochains agreeing on a subset have coboundaries agreeing there.

Note that $K_i\subseteq V_i$: the closure of an open ball of radius $\gamma_i$ lies in the open ball of radius $\gamma_0$ (resp.\ $\gamma_{i-1}$) with the same center, since $\gamma_i<\gamma_{i-1}\le\gamma_0$.

We argue by induction on $i$. For $i=0$, $t^{(0)}=\mathcal L(\eta)^{(0)}$ is defined on $U_0=\tilde U_0^{\gamma_0}$, and $d^{n-1}t^{(0)}=\eta|_{U_0}$ by $(\ast)$. For inductive step, let $i\ge1$ and assume the claim for $i-1$.

\emph{Step 1 (the difference is a cocycle).} By the induction hypothesis $t^{(i-1)}$ is defined on $\tilde U_{i-1}^{\gamma_{i-1}}$ with $d^{n-1}t^{(i-1)}=\eta$ there, and by $(\ast)$, $\eta^{(i)}=\mathcal L(\eta)^{(i)}$ is defined on $U_i$ with $d^{n-1}\eta^{(i)}=\eta$ there. Hence $\tilde t^{(i)}$ is defined on $V_i$ and $d^{n-1}\tilde t^{(i)}=\eta-\eta=0$, i.e.\ $\tilde t^{(i)}\in(\ker d^{n-1})(V_i)$.

\emph{Step 2 (local lift).} Since $V_i\subseteq U_i$, $(\ast)$ gives $\theta:=\mathcal L(\tilde t^{(i)})^{(i)}\in\check C^{n-2}(\mathcal U,\mathcal F_{EG})(V_i)$ with $d^{n-2}\theta=\tilde t^{(i)}$.

\emph{Step 3 (extension).} Let $\tilde\theta:=\mathbf{Ext}^{U_i^{\gamma_0}}_{V_i}\big(\theta;(g_i,K_i)\big)$, i.e.\ $\tilde\theta=h(\theta,g_i)$ on $V_i$ and $\tilde\theta=e_{EG}$ on $U_i^{\gamma_0}\setminus V_i$. Let $S:=\operatorname{supp}(g_i)\cap U_i^{\gamma_0}$, which is closed in $U_i^{\gamma_0}$. By (B2), $U_i^{\gamma_0}=V_i\cup(U_i^{\gamma_0}\setminus S)$ is a union of two open sets. On $V_i$, $\tilde\theta=h(\theta,g_i)$ is continuous. On $V_i\setminus S$ we have $g_i=0$, so $h(\theta,0)=e_{EG}$ (Proposition~\ref{Contractibility of EG}). Thus $\tilde\theta\equiv e_{EG}$ on $U_i^{\gamma_0}\setminus S$, and $\tilde\theta$ is continuous. By (B1), $\tilde\theta=h(\theta,1)=\theta$ on $K_i$.

Set $\bar t^{(i)}:=d^{n-2}\tilde\theta$. Then $d^{n-1}\bar t^{(i)}=0$ on $U_i^{\gamma_0}$, and by pointwise computation of $d$, $\bar t^{(i)}=d^{n-2}\theta=\tilde t^{(i)}\quad\text{on }K_i.$

\emph{Step 4 (gluing).} On $U_i^{\gamma_i}$ we have $d^{n-1}\big(\eta^{(i)}+\bar t^{(i)}\big)=\eta+0=\eta$. On the overlap $U_i^{\gamma_i}\cap\tilde U_{i-1}^{\gamma_i}\subseteq K_i$, $\eta^{(i)}+\bar t^{(i)}=\eta^{(i)}+\tilde t^{(i)}=t^{(i-1)}.$
Since $\tilde U_{i-1}^{\gamma_i}\subseteq\tilde U_{i-1}^{\gamma_{i-1}}$, the restriction $t^{(i-1)}|_{\tilde U_{i-1}^{\gamma_i}}$ is defined and satisfies $d^{n-1}t^{(i-1)}=\eta$ there. The gluing axiom thus yields $t^{(i)}$ on $\tilde U_{i-1}^{\gamma_i}\cup U_i^{\gamma_i}=\tilde U_i^{\gamma_i}$ with $d^{n-1}t^{(i)}=\eta$, completing the induction.

Finally, $\tilde U_{N-1}^{\gamma_{N-1}}=B$ because $\mathcal U^{\gamma_{N-1}}$ covers $B$.
\end{proof}

\subsection{The Eilenberg-MacLane Algorithm}\label{sec: Eilenberg-MacLane Algorithm}

Having established, in Sections~\ref{sec:bockstein-lifting}--\ref{sec:from-simplicial-to-cech}, how a persistent cohomology class $[\boldsymbol{\eta}]\in PH^n(\mathcal R(X);\mathbb Z_q)$ is transported into a cocycle $\eta\in\ker d^n_B\subseteq\check C^n(\mathcal U,\mathcal F_G)(B)$, and, in Section~\ref{sec: the-Q^n-Algorithm}, how $Q^n$ produces from such a cocycle an explicit preimage under $d^{n-1}$, we are now in a position to assemble the complete Eilenberg-MacLane Algorithm: starting from the original data $X$ and a chosen persistence interval, one proceeds through Bockstein lifting, the passage from simplicial to \Cech cohomology, and then repeated calls to $Q^n,Q^{n-1},\ldots,Q^1$ along the tower $G,BG,\ldots,B^{n-1}G$ — each call's output passed, via $i_\#$ and $j_\#$, as the input to the next — to produce the coordinate function $f_\eta: X^{(\alpha)}\to B^nG$.

\begin{algorithm}[htbp] \label{alg: Eilenberg-MacLane Algorithm}
\caption{The Eilenberg-MacLane Algorithm: $PH^n(\mathcal{R}(X); \mathbb{Z}_q)\to [X^{(\alpha)}, B^nG]$}
\begin{algorithmic}[1]
\Require Data $X$, a persistent cohomology class $[\boldsymbol{\eta}] \in PH^n(\mathcal{R}(X); \mathbb{Z}_q)$ with persistence interval $[a,b)$ satisfying $2a<b$, a scale $\alpha\in[a,b/2)$
\Require A coefficient group $G$, suggested by Bockstein lifting (Section~\ref{sec:bockstein-lifting}) applied to the representative $[\eta_G]\in H^n(\mathcal R_{2\alpha}(X);G)$
\State Define $[\eta_{\mathcal{F}_G}]\in \check{H}^n(\mathcal{U}, \mathcal{F}_G)$ by $\eta_{i_0...i_n}=\text{constant function at }\eta_G(\{x_{i_0}, ..., x_{i_n}\}) \text{ over } U_{i_0...i_n}$.
\State Set $\eta^{(n)}:=\eta_{\mathcal{F}_G}$.
\For{$k=n, n-1, ..., 1$}
    \State $\eta^{(k-1)}:= j_{\#}\Bigl(Q^{k}\bigl(i_{\#}(\eta^{(k)})\bigr)\Bigr)$
\EndFor
\State \Return $f_\eta: X^{(\alpha)}\to B^nG$ by

$\forall b\in X^{(\alpha)}$, find $j$ such that $d(b, x_j)<\alpha$. Then define $f_\eta(b):=\eta^{(0)}_j(b)\in B^nG$.
\end{algorithmic}
\end{algorithm}

The algorithm above calls $Q^n,Q^{n-1},\ldots,Q^1$ in succession, reaching $\eta^{(1)}\in\check C^1(\mathcal U,\mathcal F_{B^{n-1}G})(B)$ after the second-to-last iteration and terminating with $Q^1$ to produce $f_\eta$ directly. As discussed in Section~\ref{sec: Milnor-variant}, this last step may be replaced: rather than invoking $Q^1$, one may instead treat $\eta^{(1)}$ as a \Cech 1-cocycle representing a principal $B^{n-1}G$-bundle, and pass to its classifying map in the Milnor model $\mathcal B(B^{n-1}G)$ via the correspondence \eqref{eq:milnor-correspondence}. We record this variant below.

\begin{algorithm}[htbp] \label{alg: Milnor-Variant}
\caption{The Milnor Variant of Eilenberg-MacLane Algorithm: $PH^n(\mathcal{R}(X); \mathbb{Z}_q)\to [X^{(\alpha)}, \mathcal{B}B^{n-1}G]$}
\begin{algorithmic}[1]
\Require Data $X$, a persistent cohomology class $[\boldsymbol{\eta}] \in PH^n(\mathcal{R}(X); \mathbb{Z}_q)$ with persistence interval $[a,b)$ satisfying $2a<b$, a scale $\alpha\in[a,b/2)$
\Require A coefficient group $G$, suggested by Bockstein lifting (Section~\ref{sec:bockstein-lifting}) applied to the representative $[\eta_G]\in H^n(\mathcal R_{2\alpha}(X);G)$
\Require A partition of unity $\{\phi_x: X^{(\alpha)}\to [0,1]\}_{x\in X}$ of $X^{(\alpha)}$ subordinated to $\{B_\alpha(x)\}_{x\in X}$. Suggested: $\phi_x(y):= \frac{|\alpha - d(x,y)|_+}{\sum_{x'\in X} |\alpha - d(x',y)|_+}$, where $|\lambda|_+:=\max\{\lambda, 0\}$.
\State Define $[\eta_{\mathcal{F}_G}]\in \check{H}^n(\mathcal{U}, \mathcal{F}_G)$ by $\eta_{i_0...i_n}=\text{constant function at }\eta_G(\{x_{i_0}, ..., x_{i_n}\}) \text{ over } U_{i_0...i_n}$.
\State Set $\eta^{(n)}:=\eta_{\mathcal{F}_G}$.
\For{$k=n, n-1, ..., 2$}
    \State $\eta^{(k-1)}:= j_{\#}\Bigl(Q^{k}\bigl(i_{\#}(\eta^{(k)})\bigr)\Bigr)$
\EndFor
\State \Return Define $f_\eta: X^{(\alpha)}\to \mathcal{B}(B^{n-1}G)$ by

$\forall b\in X^{(\alpha)}$, define $\sigma_b:=\{k\ |\ d(b, x_k)<\alpha\}=\{\sigma_0,\sigma_1,...,\sigma_l\}$; Select any $j\in \sigma_b$; then 
$f_\eta(b) := \bigl( (\eta^{(1)}_{j\sigma_0}(b), \eta^{(1)}_{j\sigma_1}(b), ..., \eta^{(1)}_{j\sigma_l}(b)), (\phi_{\sigma_0}(b), \phi_{\sigma_1}(b), ..., \phi_{\sigma_l}(b))\bigr)\in \mathcal{B}(B^{n-1}G)$.
\end{algorithmic}
\end{algorithm}

An implementation of both the canonical algorithm and the Milnor variant is publicly available; see the Code Availability section at the end of this paper.

\section{Experiments}

This section presents experiments illustrating the Eilenberg-MacLane Algorithm and its Milnor variant on different examples of spaces $\mathbb X$:

\begin{enumerate}
    \item[1.] $\mathbb X=S^1$, generating both a $B\mathbb Z$-coordinate and a $B\mathbb Z_q$-coordinate, illustrating the role of Bockstein lifting in the choice of coefficient group;
    \item[2.] $\mathbb X=S^2$, generating both a $B^2\mathbb Z$-coordinate (via the canonical algorithm) and a $\mathcal B(B\mathbb Z)$-coordinate (via the Milnor variant), from the same underlying class;
    \item[3.] $\mathbb X=M(\mathbb Z_q,1)$ (a Moore space), generating a $\mathcal B(B\mathbb Z_q)$-coordinate — a new coordinate type not addressed in any prior work.
\end{enumerate}

In each experiment, we sample $\mathbb X$ to obtain a finite dataset $X$, fix an open cover of $X^{(\alpha)}$ together with the associated filtration, and apply either the Eilenberg-MacLane Algorithm or its Milnor variant to compute the resulting coordinate. This yields, for every sample point in $X$, a corresponding element of $B^nG$ or $\mathcal B(B^{n-1}G)$. In principle, a coordinate can be evaluated at any query point of $X^{(\alpha)}$, not merely at the sampled points themselves; since all experiments below use uniform sampling of $\mathbb X$, however, there is no need to query $X^{(\alpha)}$ at points other than those already sampled. The methods of this paper nonetheless permit such queries whenever desired.\\

\textit{1. Circular Coordinates on the Noisy Circle} ($\mathbb{X}=S^1$, asking $B\mathbb{Z}$/$\mathcal{B}\mathbb{Z}$-coordinates)

We begin by examining the basic properties of the Eilenberg-MacLane Algorithm, comparing it against the classical, principal-bundle-based degree-one algorithm underlying circular coordinates \cite{desilva2011persistent,perea2019sparse}. The experiment starts with a canonical example: the noisy circle. We sample $1000$ points with noise from $\mathbb X=S^1$ to form the dataset $X$, and fix the simplest possible open cover consisting of three balls, $\mathcal U=\{B_\alpha(v_0),B_\alpha(v_1),B_\alpha(v_2)\}$, where distances are measured by the geodesic distance on the circle. The centers $v_0,v_1,v_2$ are the three cube roots of unity, at $0^\circ,120^\circ,240^\circ$ respectively, i.e. $v_0=(1,0)$, $v_1=(-\tfrac12,\tfrac{\sqrt3}2)$, $v_2=(-\tfrac12,-\tfrac{\sqrt3}2)$, and $\alpha:=d(v_0,v_1)\approx2.0944$ is the biggest radius for which the three balls form a good cover on $S^1$.

\begin{figure}[htbp]
    \centering
    \includegraphics[width=0.32\textwidth]{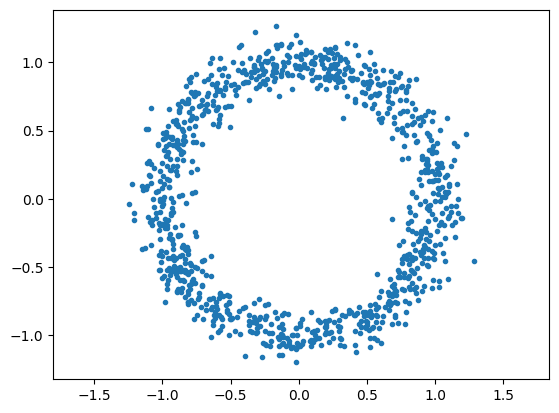}
    \includegraphics[width=0.31\textwidth]{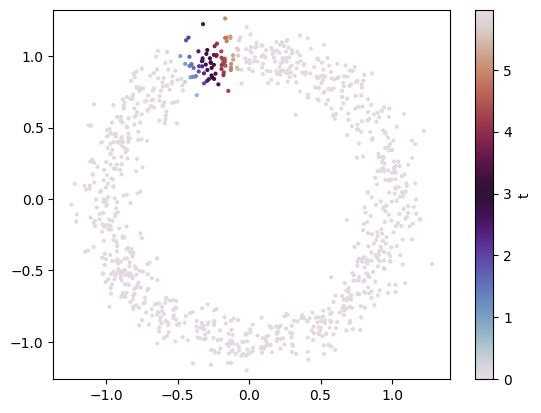}\includegraphics[width=0.31\textwidth]{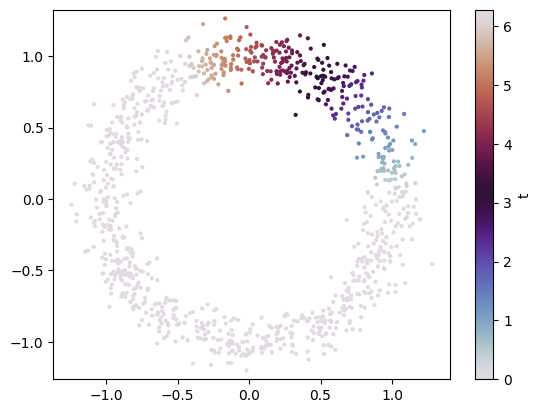}
    \caption{Comparison of the original sampling from $S^1$ (left), the canonical Eilenberg-MacLane Algorithm (middle), and the Milnor variant (right). The colormap assigns to each data point $x\in X$ a color in $[0,2\pi)$, representing $f_\eta(x)\in B\mathbb Z\simeq\mathcal B\mathbb Z\simeq S^1$.}
    \label{fig:exp1-s1-comparison}
\end{figure}

We use the radii-decreasing sequence $\gamma=(\gamma_0,\gamma_1,\gamma_2)=(\alpha,1.6,1.2)$ explained in Section~\ref{sec: the-Q^n-Algorithm}. Beyond the mathematical requirement that $\mathcal U^{\gamma_2}$ still cover $S^1$, this choice was also made for practical reasons for which each pairwise intersection $U_i\cap U_j$ retains enough sampled points that the resulting figures remain legible. Under this setup, we generate Eilenberg-MacLane coordinates — both the canonical algorithm and the Milnor variant — from
$$\eta:=\{\eta_{01}\equiv1,\ \eta_{12}=\eta_{02}\equiv0\}\in\check Z^1(\mathcal U,\mathcal F_{\mathbb Z}),$$
the standard generator of $H^1(S^1;\mathbb Z)\cong\mathbb Z$. Since $n=1$ here, the Milnor variant coincides with the classical circular coordinate construction. Figure~\ref{fig:exp1-s1-comparison} displays the original dataset alongside the canonical Eilenberg-MacLane coordinate and the Milnor variant.

\begin{figure}[htbp]
    \centering
    \includegraphics[width=0.6\textwidth]{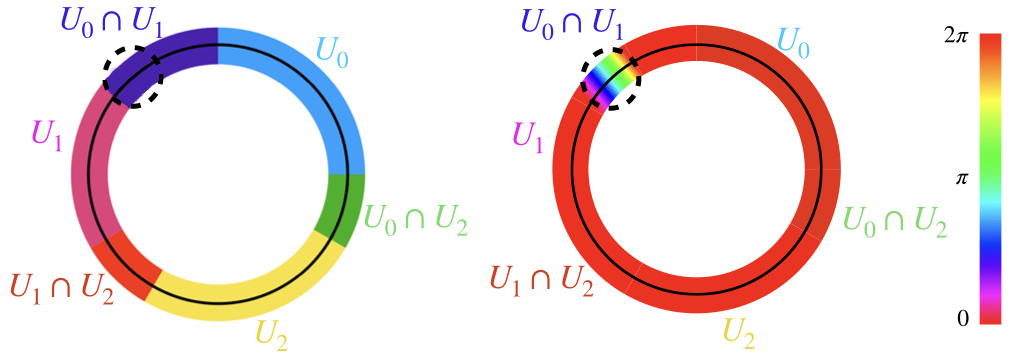}
    \caption{The region where the canonical Eilenberg-MacLane algorithm produces nontrivial bars. As indicated by the dashed circle, it locates at $U_0$ between radius $\gamma_0$ and $\gamma_1$, on the $U_1$ side.}
    \label{fig: S1-nontrivial-bars-representation}
\end{figure}

As we saw in Section~\ref{sec: Milnor-variant}, the Milnor variant requires the construction of a partition of unity, a more demanding condition than the bump functions $g_i$ required by the canonical Eilenberg-MacLane algorithm; the canonical line is therefore evidently more broadly applicable. A notable difference, however, is where the nontrivial bars of the Eilenberg-MacLane coordinate occur: in the canonical line, they occur on the portion added when $\tilde U_{i-1}^{\gamma_{i-1}}$ is extended to $\tilde U_i^{\gamma_i}$. More precisely, they occur on the annulus between the two radii, $\gamma_i$ and $\gamma_{i-1}$, of $\tilde U_{i-1}$, intersected with $U_i^{\gamma_0}$, the domain of the newly introduced $\eta^{(i)}$. Figure~\ref{fig: S1-nontrivial-bars-representation} shows the region where the canonical Eilenberg-MacLane algorithm produces genuinely nontrivial bars. By contrast, since the Milnor variant uses a partition of unity, the region where nontrivial bars occur is the support of some $\phi_x$, i.e., some $U_i$.

Both the canonical line and the Milnor variant exhibit a similar \textit{geometric distortion}, which is suggested to be improved by algorithms such as harmonic smoothing, as indicated in \cite{desilva2011persistent,perea2019sparse}. Beyond this, we will also see in the experiments below that a geometric distortion of this kind can be corrected for using methods such as Isomap and MDS.

\textit{2. Bockstein Lifting on the Circle} ($\mathbb{X}=S^1$, asking $B\mathbb{Z}$/$B\mathbb{Z}_q$-coordinates)

\begin{figure}[htbp]
    \centering
    \includegraphics[width=0.32\textwidth]{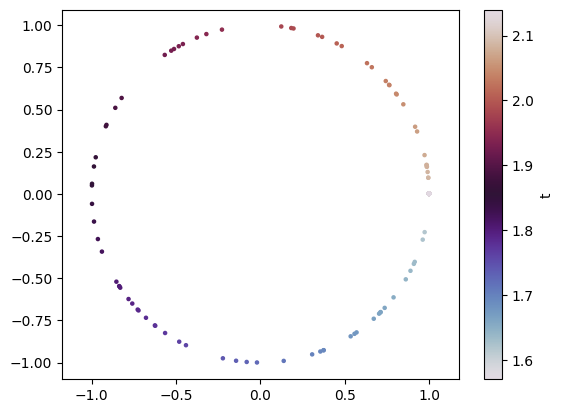}
    \includegraphics[width=0.32\textwidth]{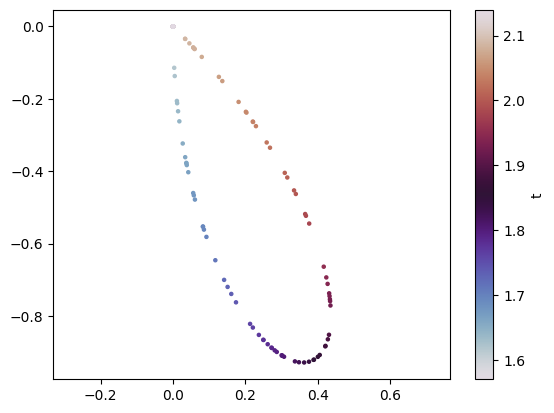}
    \includegraphics[width=0.32\textwidth]{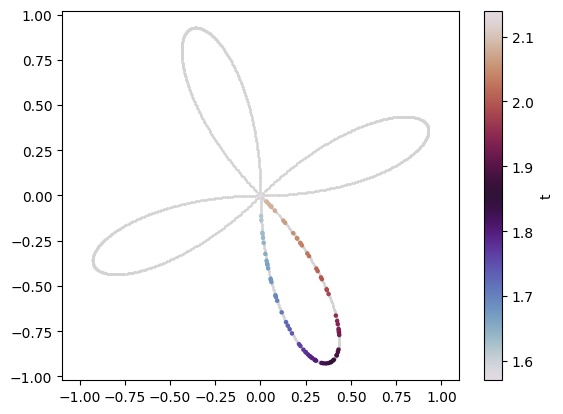}
    \caption{The same dataset $X$ (noisy circle) represented under the $B\mathbb Z$-coordinate (left) and the $B\mathbb Z_q$-coordinate with $q=5$ (middle). The $B\mathbb Z$-coordinate of $X$ is embedded in $S^1$, while the $B\mathbb Z_q$-coordinate of $X$ is embedded in the $1$-skeleton of $L_q^\infty$, namely $\bigvee_{i=0}^{q-1}S^1$. The right panel shows the same embedding, with the entire codomain $\bigvee_{i=0}^{q-1}S^1$ additionally displayed in gray.}
    \label{fig: S1-BZ_q-coordinate}
\end{figure}

We next illustrate the role of Bockstein lifting. For the same dataset $X$ sampled from $S^1$, we generate both a $B\mathbb Z$-coordinate and a $B\mathbb Z_q$-coordinate (with $q=5$), and plot both alongside the dataset in Figure~\ref{fig: S1-BZ_q-coordinate}. The $B\mathbb Z$-coordinate is converted, via a homotopy equivalence, into an element of $S^1$. The $B\mathbb Z_q$-coordinate (lens coordinate), since the nontrivial bars it produces never exceed degree $1$, lies entirely within the $1$-skeleton of $K(\mathbb Z_q,1)=L_q^\infty$, namely $\bigvee_{i=0}^{q-1}S^1$. Although the image of the $B\mathbb Z_q$-coordinate of $X$ occupies only a single $S^1$ summand of $L_q^\infty$, a consequence of $X$ itself having been sampled from $S^1$, this example illustrates that different coordinates need not carry the same interpretability or quality of representation. Here, Bockstein lifting identifies the underlying coefficient group of $\eta$ in advance, leading one to prefer the $B\mathbb Z$-coordinate — the more informative of the two — rather than generating a coordinate valued in the more complex space $L_q^\infty$.\\

\textit{3. Validation on Iterated Bars} ($\mathbb X=S^2$, asking $B^2\mathbb Z$-coordinates)

We next validate the Eilenberg-MacLane Algorithm (Algorithm~\ref{alg: Eilenberg-MacLane Algorithm}), and in particular its correctness in successively lowering the degree of \Cech cohomology from dimension $2$ down to dimension $0$. Motivated by the geometric distortion illustrated in Experiment 1, we sample $4000$ points from the sphere $S^2$, while additionally sampling $800$ points as a local density concentrated on the specific regions where the bump functions $g_i$ are nontrivial. We fix the simplest open cover consisting of four balls, $\mathcal U=\{B_\alpha(v_0),\ldots,B_\alpha(v_3)\}$, where distances are measured by the geodesic distance on the sphere. The centers $v_0,v_1,v_2,v_3$ are the four vertices of a regular tetrahedron inscribed in the sphere, i.e., $(0,0,1)$, $(\tfrac{2\sqrt2}3,0,-\tfrac{1}3)$, $(-\tfrac{\sqrt2}3,\tfrac{\sqrt6}3,-\tfrac{1}3)$, $(-\tfrac{\sqrt2}3,-\tfrac{\sqrt6}3,-\tfrac{1}3)$. We take $\alpha:=d(v_0,v_1)\approx1.9106$ as the largest radius for which the four balls form a good cover of $S^2$.

We use the radii-decreasing sequence $\gamma=(\gamma_0,\gamma_1,\gamma_2,\gamma_3)=(\alpha,1.8,1.7,1.6)$ as explained in Section~\ref{sec: the-Q^n-Algorithm}, so that $\mathcal U^{\gamma_3}$ still covers $S^2$. Under this setup, we generate Eilenberg-MacLane coordinates from
$$\eta:=\{\eta_{012}\equiv1,\ \eta_{013}=\eta_{023}=\eta_{123}\equiv0\}\in\check Z^2(\mathcal U,\mathcal F_{\mathbb Z}),$$
the standard generator of $H^2(S^2;\mathbb Z)\cong\mathbb Z$.

For the resulting $4800$ coordinates in $B^2\mathbb Z$, we convert them into elements of $B^2\mathbb{Z} = B(B\mathbb Z)\simeq BS^1$ via the homotopy equivalence
$$E\mathbb Z\simeq\mathbb R;\qquad [(g_0,\ldots,g_n),(t_0,\ldots,t_n)]\ \mapsto\ \sum_{k=0}^n(t_0+\cdots+t_k)g_k,$$
obtaining $4800$ bars of the form $[(z_0,z_1),(t_0,t_1,t_2)]$. (Recall that, in the bar construction of Definition~\ref{Bar Construction}, the action of $G$ on $EG$ is defined by right multiplication on the last coordinate $g_n$; consequently, an element of $BG$ may be represented as $[(g_0,\ldots,g_{n-1}),(t_0,\ldots,t_{n-1},t_n)]$, with the last coordinate $g_n$ omitted, since it carries no information in $BG$.) Since the dataset is sampled from $S^2$, we observe that $z_0=e^{2\pi i(t_0+t_1)}$ and $z_1=e^{2\pi i(1-t_0)}$; consequently, each transformed coordinate landing in $BS^1$ is uniquely represented by the pair $(t_0,t_1)$.

\begin{figure}[htbp]
    \centering
    \includegraphics[width=0.33\textwidth]{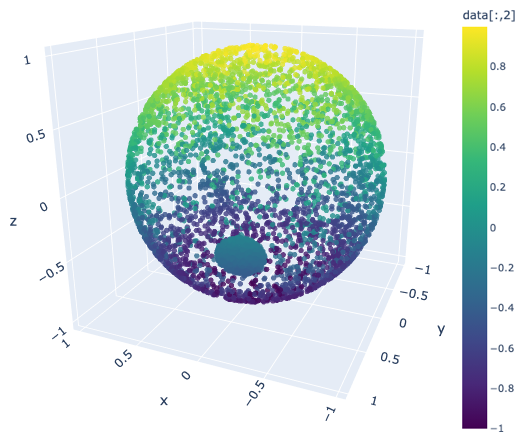}
    \includegraphics[width=0.35\textwidth]{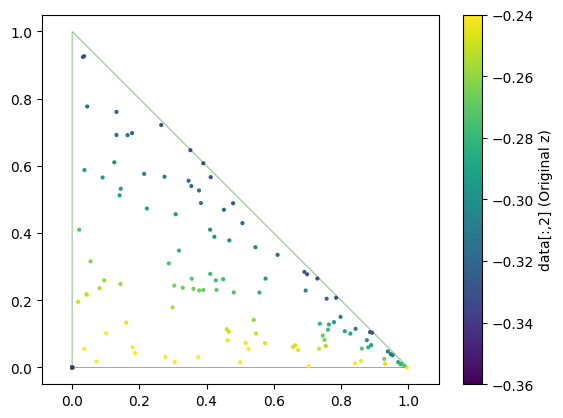}\includegraphics[width=0.3\textwidth]{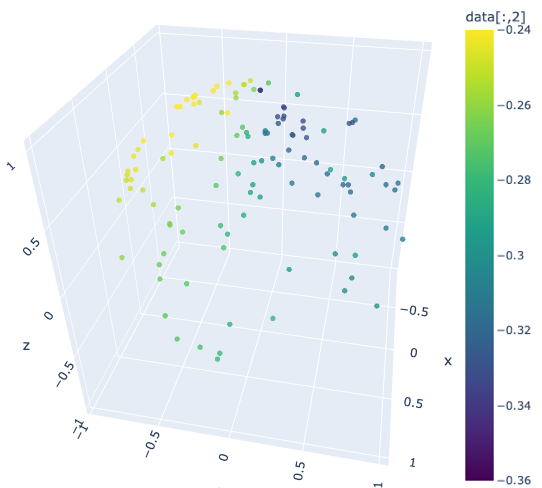}
    \caption{The dataset $S^2$ (left), and the corresponding point cloud of Eilenberg-MacLane coordinates $[(z_0,z_1),(t_0,t_1,t_2)]\in BS^1\simeq B^2\mathbb Z$, plotted via $(t_0,t_1)$ (middle), colored by the original $z$. Since the equivalence relation defining the bar construction identifies every point on the green triangle as a single point, the right panel shows the result of homeomorphically mapping this point cloud of coordinates onto $S^2$ via stereographic projection.}
    \label{fig:S2_BBZ-coordinates}
\end{figure}

Figure~\ref{fig:S2_BBZ-coordinates} displays the original dataset $S^2$ together with the resulting pairs $(t_0,t_1)$ obtained from the Eilenberg-MacLane coordinate. Note that, since $[(z_0,z_1),(t_0,t_1,t_2)]\in BS^1$, the equivalence relation defining the bar construction identifies every point on the green triangle in the figure as the same point; we may therefore further map this space homeomorphically to $S^2$ via stereographic projection, yielding the third panel of Figure~\ref{fig:S2_BBZ-coordinates}. This panel shows the new point cloud obtained in coordinate space from the original point cloud, with a colormap on the original $z$ illustrating the correspondence between the original and new point clouds induced by $f_\eta$.

This experiment demonstrates that our proposed pipeline can iteratively lower the degree of cohomology. Although the results show that the geometric distortion exhibited in Experiment 1 becomes more pronounced under iteration, this issue can be partially corrected using the Milnor-variant-based pipeline introduced next.\\

\textit{4. Recovering the Distortion via the Milnor Variant and Isomap} ($\mathbb X=S^2$, asking $\mathcal{B}(B\mathbb Z)$-coordinates)

One can partially recover the geometric distortion illustrated in Experiment 1 via the Milnor variant and Isomap. Without the local density sampling used in Experiment 3, we in this experiment compute a $\mathcal B(B\mathbb Z)$-coordinate on the same $4000$ points sampled from $S^2$, using the Milnor variant pipeline, and equip it with a metric.

\begin{figure}[htbp]
    \centering
    \includegraphics[width=0.31\textwidth]{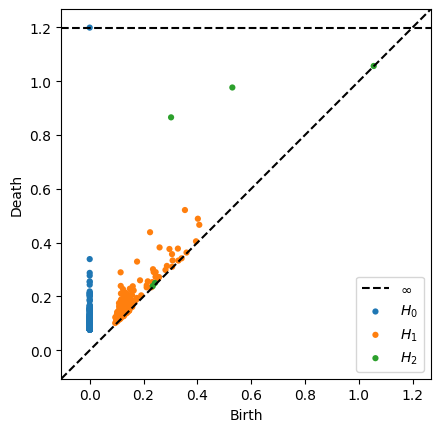}
    \includegraphics[width=0.37\textwidth]{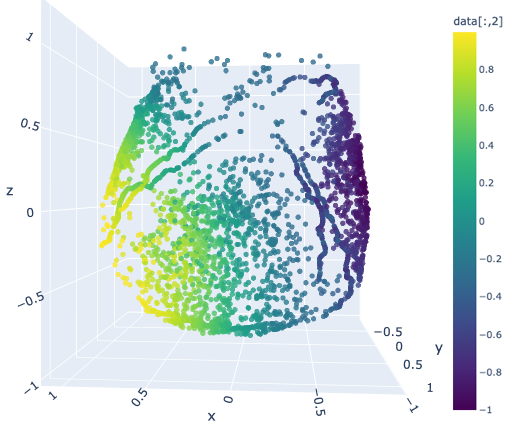}\includegraphics[width=0.31\textwidth]{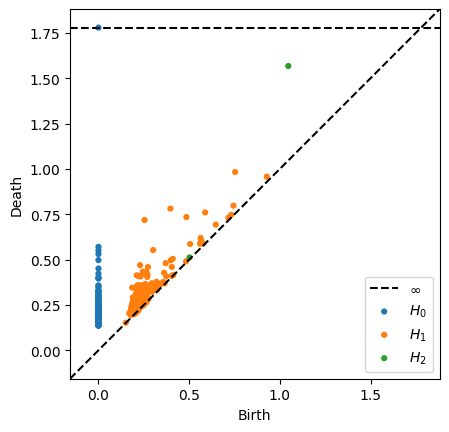}
    \caption{Left: the persistence diagram computed on $(\mathbb{CP}^\infty,d_{\mathbb{CP}^\infty})\simeq\mathcal B(B\mathbb Z)$, obtained by applying the Milnor variant to $\eta\in \check{H}^2(\mathcal{U};\mathcal{F}_\mathbb{Z})$ on a $4000$-point sampling from $S^2$. After obtaining these coordinates via the Milnor variant, we run Isomap on $(\mathbb{CP}^\infty,d_{\mathbb{CP}^\infty})$ and retain the first three dimensions, yielding the persistence diagram on the right; the resulting point cloud (colored by original $z$), embedded in $\mathbb R^3$ via Isomap, is also shown, successfully recovering the geometry of the original point cloud $S^2$.}
    \label{fig: S^2-Milnor-Variant}
\end{figure}

With the translation $\mathcal B(B\mathbb Z)\simeq\mathcal BS^1\simeq\mathbb{CP}^\infty$, the second step is given by
$$\mathcal ES^1\simeq S^\infty;\qquad [(g_0,g_1,\ldots,g_n),(t_0,t_1,\ldots,t_n)]\mapsto(\sqrt{t_0}g_0,\sqrt{t_1}g_1,\ldots,\sqrt{t_n}g_n,0,0,\ldots).$$
This yields the coordinates of the $4000$ sample points as points in $\mathbb{CP}^\infty$. We then equip $\mathbb{CP}^\infty$ with the Fubini-Study metric:
$$d_{\mathbb{CP}^\infty}([z_0:z_1:\cdots:z_n],[w_0:w_1:\cdots:w_n]):=\arccos\left(\left|\sum_{i=0}^nz_i\overline{w_i}\right|\right).$$
To understand the resulting point cloud of coordinates, we compute its persistence diagram (via Ripser) with respect to $(\mathbb{CP}^\infty,d_{\mathbb{CP}^\infty})$. The result is shown in Figure~\ref{fig: S^2-Milnor-Variant}. We observe two dominant persistence classes in dimension $2$ — rather than the single class expected from $H^2(S^2;\mathbb Z)\cong\mathbb Z$ — indicating that the coordinates produced by the Milnor variant pipeline fail to faithfully recover the metric, resulting in a distortion of the underlying topology.

We therefore seek to embed the coordinates lying in $(\mathbb{CP}^\infty,d_{\mathbb{CP}^\infty})$ into some finite-dimensional Euclidean space. The standard approach for this is Isomap. Applying Isomap to the point cloud of coordinates (with, e.g., $10$ neighbors and $8$ components), we obtain an explained variance, for the first eight dimensions, of $[0.342,0.642,0.893,0.927,0.954,0.971,0.986,1.000]$, indicating that the data can be reduced to approximately its first three dimensions, which together retain $89\%$ of the variance. The middle panel of Figure~\ref{fig: S^2-Milnor-Variant} shows the resulting point cloud after embedding into $3$-dimensional Euclidean space via Isomap, in which nearly all of the geometry of $S^2$ is preserved. Applying Ripser again to the new distance matrix obtained after Isomap, the persistence diagram is fully restored to the expected geometric structure.

\textit{5. A Genuinely New Coordinate Type: $\mathcal B(B\mathbb Z_q)$-coordinates on $\mathbb{RP}^2$} ($\mathbb X=M(\mathbb{Z}_q,1 )$, asking $\mathcal{B}(B\mathbb Z_q)$-coordinates)

Experiment 4 distilled a mature pipeline for visualizing the coordinate space concretely: Starting from a dataset $X$, we first use the Eilenberg-MacLane Algorithm (the canonical pipeline or the Milnor variant) to generate a coordinate $f_\eta(X)\subseteq B^nG$ or $\mathcal B(B^{n-1}G)$, and then equip this coordinate space with a metric. Next, via Isomap we obtain a new point cloud embedded in $\mathbb R^N$, which can be visualized or topologically verified by computing its persistence diagram.

We now apply this pipeline to generate $\mathcal B(B\mathbb Z_q)$-coordinates (taking $q=2$ as our example) — a coordinate type that has not, to our knowledge, appeared in any prior work. It is worth emphasizing that this coordinate is not a lens coordinate of the kind already treated by Polanco and Perea~\cite{polanco2019coordinatizingdatalensspaces}, who produce $K(\mathbb Z_q,1)$-valued (lens) coordinates from a class in $H^1(-;\mathbb Z_q)$ — for instance, on $\mathbb{RP}^2\cong M(\mathbb Z_2,1)$ itself, whose $H^1(\mathbb{RP}^2;\mathbb Z_2)\cong\mathbb Z_2$ falls squarely within their framework, nor is it a real or complex projective coordinate as in Perea~\cite{perea2017multiscale}. What we construct instead is a genuinely \emph{degree-two} coordinate, valued in $K(\mathbb Z_q,2)=\mathcal B(B\mathbb Z_q)$, obtained from the generator of $H^2(M(\mathbb Z_q,1);\mathbb Z_q)$.

Recall that the cohomology groups of $M(\mathbb Z_q,1)$ (Moore space) are as follows:

$$H^i(M(\mathbb Z_q,1);\mathbb Z) = \left\{\begin{array}{ll}
    \mathbb{Z} & (i=0) \\
    0 & (i=1)\\
    \mathbb{Z}_q & (i=2) \\
    0 & (i\geq 3)
\end{array}\right. \qquad
H^i(M(\mathbb Z_q,1);\mathbb Z_q) = \left\{\begin{array}{ll}
    \mathbb{Z}_q & (i=0) \\
    \mathbb{Z}_q & (i=1)\\
    \mathbb{Z}_q & (i=2) \\
    0 & (i\geq 3)
\end{array}\right.
$$

That is, the $\mathbb Z_q$ appearing in degree $2$ is native to $M(\mathbb Z_q,1)$: $H^2(M(\mathbb Z_q,1);\mathbb Z)\cong\mathbb Z_q$ is itself a $q$-torsion group, so the class detected in $\mathbb Z_q$ coefficients genuinely originates from this torsion, rather than being a mod-$q$ artifact of some richer integral structure. This provides ample motivation to work directly with the $\mathbb Z_q$-coordinate at this degree. We accordingly use $M(\mathbb Z_q,1)$ as our test space and ask for its $\mathcal B(B\mathbb Z_q)$-coordinate. In practice we take $q=2$, so $M(\mathbb Z_q,1)=\mathbb{RP}^2$.

In this experiment, we choose the typical $6$-vertex minimal triangulation of $\mathbb{RP}^2$, which consists of $10$ faces, $15$ edges, and $6$ vertices, as shown in the left panel. We uniformly sample $4000$ points from this triangulation, encoding each point via barycentric coordinates $((v_{i_0},v_{i_1},v_{i_2}),(t_{i_0},t_{i_1},t_{i_2}))$, and equip this model with the simplicial $\ell^1$-distance intrinsic to the triangulation. Concretely, given two such points, we first extend each to a weight assignment $\{\text{vertex index}\mapsto\text{weight}\}$, where vertices not appearing in the point's representative are assigned weight $0$; the distance between two points is then defined as the sum of the absolute differences of their weight assignments, i.e.
$$d_{\ell^1}(x,y):=\sum_{i=0}^5|w_x(v_i)-w_y(v_i)|.$$
Taking each vertex $v$ as the center of an open star of radius $r$, i.e. $B_\alpha(v):=\{x\mid d_{\ell^1}(x,v)<r\}$, we find that the union of such balls forms a good cover for radii ranging between $4/3$ and $2$. We accordingly choose the radii-decreasing sequence $\gamma=(\gamma_0,\ldots,\gamma_5)=(2.0,1.9,1.8,1.7,1.6,1.5)$. We then take
$$\eta:=\{\eta_{013}\equiv1,\ \eta_{i_0i_1i_2}\equiv0\text{ for all other }(i_0,i_1,i_2)\}$$
as the generator of $H^2(M(\mathbb Z_2,1);\mathbb Z_2)$, and run the Milnor variant. This yields $\mathcal B(B\mathbb Z_2)$-coordinates of the form
$$c_1=\Bigl[\bigl(e,\ e,\ e,\ ((0.4466,0.5534),(1)),\ e,\ e\bigr),\ (0.0277,\ 0.9124,\ 0.0000,\ 0.0600,\ 0.0000,\ 0.0000)\Bigr]$$
where $e=((1),)$ is the identity of $B\mathbb{Z}_2$, and the second entry records the (rounded) barycentric weights $(t_0,\ldots,t_5)$ over the six vertices, and each nonzero coordinate carries an associated element of $\mathcal B\mathbb Z_2$ recording the corresponding inner bar (e.g., $((0.4466,0.5534),(1))$ at the vertex indexed $3$, and the identity bar elsewhere).

We observe experimentally that these inner bars (such as $e$ or $((0.4466,0.5534),(1))$ in $c_1$) never exceed degree $1$, so they may unconditionally be converted into elements of $S^1$: each such bar is either the identity bar, or of the form $((t_0,t_1),(g_0))$, which, under the identification defined by the bar construction on those $t_0,t_1$, and $g_0$'s, reduces to $\exp(2\pi it_0g_0)\in S^1$. For instance, the inner bar $((0.4466,0.5534),(1))$ appearing in $c_1$ above reduces to $\exp(2\pi i\times 0.4466 \times 1)\in S^1$, and $e$ simply converts into $1\in S^1$.

Now that our coordinates have been converted into elements of $\mathcal BS^1$, we may further use the identification $\mathcal BS^1\simeq\mathbb{CP}^\infty$ from Experiment 4, together with the Fubini-Study metric equipped on $\mathbb{CP}^\infty$, to compute the pairwise distances between coordinates, and use Isomap to embed the coordinates into Euclidean space.

\begin{figure}[htbp]
    \centering
    \includegraphics[width=0.31\textwidth]{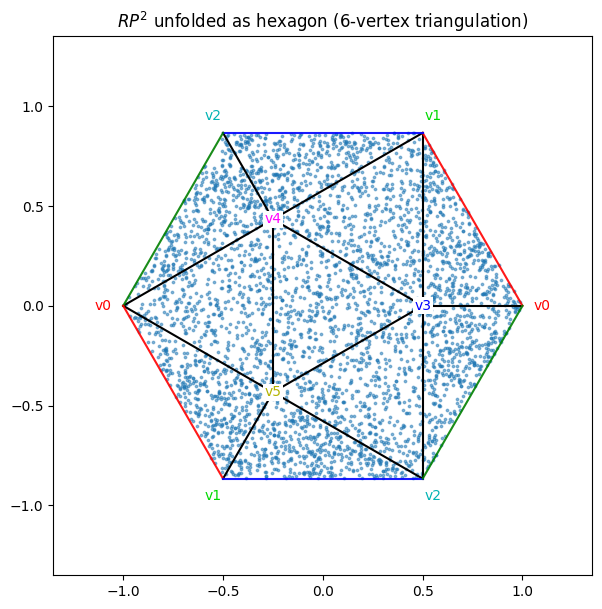}
    \includegraphics[width=0.31\textwidth]{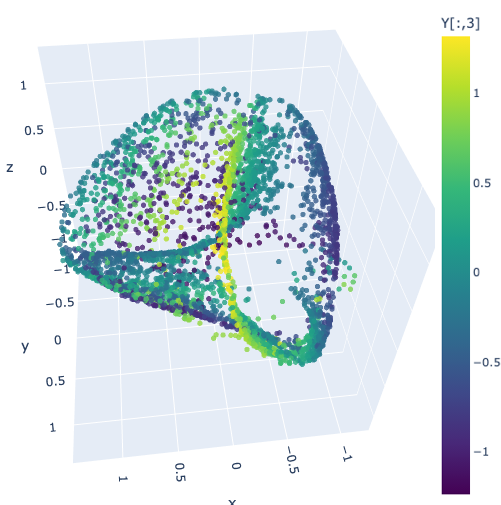}
    \includegraphics[width=0.36\textwidth]{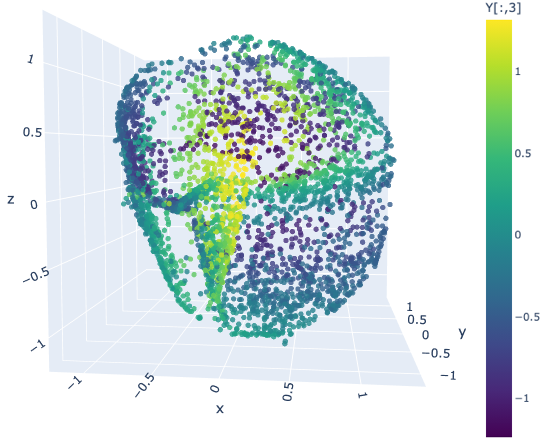}
    \\
    \includegraphics[width=0.32\textwidth]{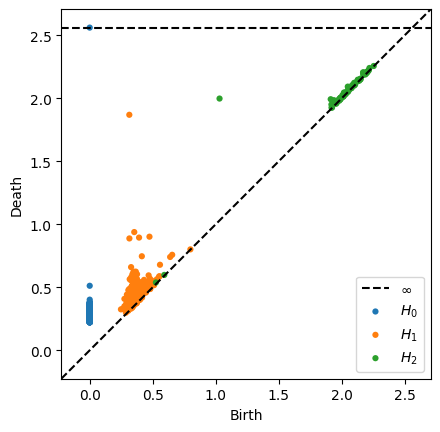}
    \includegraphics[width=0.32\textwidth]{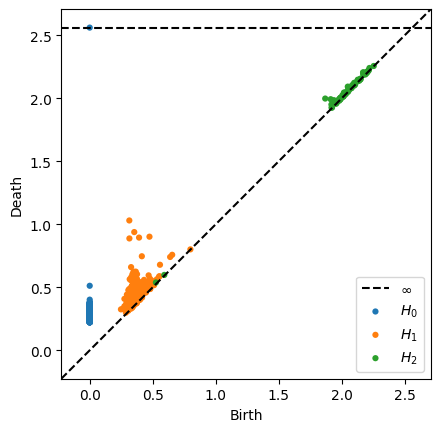}
    \caption{Top left: the triangulation of $\mathbb{RP}^2$, together with $4000$ points sampled uniformly from it, forming the point cloud $X$. Top middle and top right: the coordinates induced by the generator of $H^2(\mathbb{RP}^2;\mathbb Z_2)$ via the Milnor variant, embedded into the first three dimensions of Euclidean space after Isomap, colored by the fourth dimension. Bottom left and bottom right: the persistence diagrams obtained by applying Ripser, with coefficients in $\mathbb Z_2$ and $\mathbb Z_3$ respectively, to the new point cloud resulting from the Milnor variant followed by Isomap. Both persistence diagrams correctly confirm the cohomology of $\mathbb{RP}^2$ in the corresponding coefficients.}
    \label{fig: RP^2-Milnor-Variant}
\end{figure}

Applying Isomap to the point cloud of coordinates (with $10$ neighbors and $8$ components), we obtain an explained variance, for the first eight dimensions, of 
\[
[0.206, 0.403, 0.589, 0.773, 0.945, 0.965, 0.983, 1.000],\] indicating that the data can be reduced to approximately its first four dimensions, which together retain $77\%$ of the variance. Since these four dimensions cannot be directly visualized, the middle panel of Figure~\ref{fig: RP^2-Milnor-Variant} shows the resulting point cloud after embedding into the first $3$ of these dimensions, colored by the fourth. This is essentially the result of forcibly embedding $\mathbb{RP}^2$ into $\mathbb R^3$, known as the \href{https://en.wikipedia.org/wiki/Roman_surface}{Roman surface}. Applying Ripser again to the new distance matrix obtained after Isomap, the persistence diagram is fully restored to the expected geometric structure of $\mathbb{RP}^2$.

\section{Future Work}
\label{sec:future-work}

We conclude by outlining several directions left open by the present work.

\textbf{1. Comparing the canonical pipeline and the Milnor variant.} Both pipelines identify $\check H^1(B;\mathcal F_H)$, with $H=B^{n-1}G$, with homotopy classes of maps into a classifying space of $H$: the canonical pipeline through the connecting homomorphism of Section~\ref{sec:chains-end}, and the Milnor variant through the classical correspondence~\eqref{eq:milnor-correspondence}. It would be worthwhile to determine whether these two identifications agree (possibly up to sign) under a homotopy equivalence $BH\simeq\mathcal BH$, and to compare the resulting coordinates geometrically.

\textbf{2. Complexity analysis.} While we have argued informally that both the canonical pipeline and the Milnor variant admit polynomial-time implementations, a rigorous complexity analysis — bounding the cost of the $Q^n$-Algorithm's recursive gluing procedure, and comparing it against the partition-of-unity-based construction underlying the Milnor variant — has not been carried out here, and is left for future work.

\textbf{3. Combining multiple persistence classes via the Postnikov tower.} The present construction extracts a coordinate from a single persistent cohomology class at a time. A natural, and considerably more ambitious, extension is to combine several such classes — potentially in different degrees and with different coefficient groups — into a single coordinate that captures more of the homotopy type of the underlying space than any individual Eilenberg-MacLane coordinate can. Since a space is generally not homotopy equivalent to a product of Eilenberg-MacLane spaces, realizing such a combination would require understanding how the individual classes are glued together

\section*{Code Availability}
An implementation of the Sparse Eilenberg-MacLane Algorithm described in this paper, together with scripts used to generate the experiments in Section~5, is available at \url{https://github.com/tonyxiaochenxiao/Eilenberg-MacLane-Coordinates}.

\appendix
\section{The contracting homotopy of $EG$}\label{app:contracting-homotopy}

We verify that the map $h$ of~\eqref{eq:contracting-homotopy} is a well-defined continuous map $EG\times[0,1]\to EG$ with $h(\cdot,1)=\mathrm{id}$ and $h(\cdot,0)\equiv e$.

\emph{Continuity.} On each piece $G^{n+1}\times\Delta^n\times[0,1]$, the map $h$ is the product of
$$G^{n+1}\to G^{n+2},\ (g_0,\dots,g_n)\mapsto((g_0\cdots g_n)^{-1},g_0,\dots,g_n)\quad\text{and,}$$
$$\quad \Delta^n\times[0,1]\to\Delta^{n+1},\ ((t_0,\dots,t_n),t)\mapsto(1-t,tt_0,\dots,tt_n),$$
which are continuous since the multiplication and inversion of $G$ are continuous. Since $[0,1]$ is locally compact Hausdorff, the product of the quotient map $\bigsqcup_{n}G^{n+1}\times\Delta^n\to EG$ with $[0,1]$ is again a quotient map, so it suffices to check that $h$ respects the relations of Definition~\ref{Bar Construction}.

\emph{Well-definedness.} Write $\pi:=g_0\cdots g_n$, and note that $h$ inserts a new first entry, so that an entry in position $i$ moves to position $i+1$.
\begin{enumerate}
    \item If $t_0=0$, then $((g_0,\dots,g_n),(t_0,\dots,t_n))\sim((g_1,\dots,g_n),(t_1,\dots,t_n))$. The image of the left side has second barycentric coordinate $tt_0=0$, so the face relation merges its first two entries into $\pi^{-1}g_0=(g_1\cdots g_n)^{-1}$, which is exactly the image of the right side.
    \item If $t_i=0$ for some $i>0$, then $((g_0,\dots,g_n),(t_0,\dots,t_n))\sim((\dots,g_{i-1}g_i,\dots),(\dots,t_{i-1},t_{i+1},\dots))$. Both sides have the same product $\pi$, and in the image of the left side the coordinate $tt_i=0$ now sits in position $i+1>0$, so the face relation merges $g_{i-1}$ and $g_i$, yielding the image of the right side.
    \item If $g_i=e$ for some $i<n$, then $((g_0,\dots,g_n),(t_0,\dots,t_n))\sim((\dots,g_{i-1},g_{i+1},\dots),(\dots,t_{i-1},t_i+t_{i+1},\dots))$. Both sides have the same product $\pi$, and in the image of the left side the entry $e$ now sits in position $i+1<n+1$, so the degeneracy relation removes it and merges $tt_i+tt_{i+1}$, yielding the image of the right side.
\end{enumerate}

\emph{Endpoints.} For $t=1$, the first barycentric coordinate is $0$, so the face relation removes the first entry and $h(x,1)=x$. For $t=0$, the barycentric coordinates are $(1,0,\dots,0)$, so the face relations merge all entries into $\pi^{-1}g_0\cdots g_n=e$, and $h(x,0)$ is the point $((e),(1))$, the identity of $EG$.

\section*{Acknowledgments}
T.X. Xiao would like to thank Ben Knudsen, Cary Malkiewich and  Tanishq Bhatia for several fruitful conversations clarifying some of the mathematical constructions in this work.

\end{document}